\documentclass[preprint]{imsart}
\RequirePackage[OT1]{fontenc}
\RequirePackage{amsthm,amsmath}
\RequirePackage[numbers]{natbib}
\RequirePackage[colorlinks,citecolor=blue,urlcolor=blue]{hyperref}
\RequirePackage{hypernat}
\usepackage{amssymb,mathrsfs}
\usepackage{amsmath}
\usepackage{amssymb,mathrsfs}
\usepackage {graphicx}
\usepackage {comment}

\usepackage{multirow}
\usepackage{tikz}
\usetikzlibrary{calc,shapes}
\usetikzlibrary{positioning}
\usepackage{multirow}
\usepackage{tabularx}

\numberwithin{equation}{section}

\newtheorem{corollary}{Corollary}[section]
\newtheorem{lemma}{Lemma}[section]
\newtheorem{assumption}{Assumption}
\newtheorem{condition}{Condition}
\newtheorem{remark}{Remark}[section]
\newtheorem{example}{Example}
\newcommand{\Thm}[1]{\begin{theorem}\label{thm.#1}}
\newcommand{\Cor}[1]{\begin{corollary}\label{cor.#1}}
\newcommand{\Prop}[1]{\begin{proposition}\label{prop.#1}}
\newcommand{\Lm}[1]{\begin{lemma}\label{lm.#1}}
\newcommand{\Ass}[1]{\begin{assumption}\label{ass.#1}}
\newcommand{\Ex}[1]{\begin{example}\label{ex.#1}\rm}
\newcommand{\Rem}[1]{\begin{remark}\label{rem.#1}\rm}
\newcommand{\thmref}[1]{Theorem~\ref{thm.#1}}

\newcommand{\lmref}[1]{Lemma~\ref{lm.#1}}

\newcommand{\remref}[1]{Remark~\ref{rem.#1}}
\newcommand{\exref}[1]{Example~\ref{ex.#1}}
\newcommand{\Cond}[1]{\begin{condition}\label{cond.#1}\rm}

\def\qed{\hfill \mbox{\raggedright \rule{0.1in}{0.1in}}}

\def\Cov{\mathop{\rm Cov}\nolimits}
\def\Var{\mathop{\rm Var}\nolimits}%
\def\R{\mathbb R}
\def\bea{$$\begin{aligned}}
\def\eea{\end{aligned}$$}

\def\el{\mathscr{R}_{n}}

\def\implies{{\Longrightarrow}}

\def\und{\quad \mbox{and} \quad}

\def\1{\mathbf 1}
\def\ave{n^{-1}\sum_{j=1}^n}

\def\Ave{\frac{1}{n}\sum_{j=1}^n}
\def\Avi{\frac{1}{n}\sum_{i=1}^n}

\def\pave{\sum_{j=1}^n \pi_j}

\def\lam{\lambda} 

\def\maxj{\max_{1\leq j \leq n}}

\def\trace{\mathrm{trace}}
\def\bmS{\tilde{\bdm{\S}}}
\def\bms{\tilde{\bdm{\s}}}
\renewcommand{\Ref}[1]{(\ref{#1})}
\newcommand{\bel}[1]{\begin{equation}\label{#1}}
\newcommand{\tel}[1]{\begin{equation}\tag{#1}}
\newcommand{\ee}{\end{equation}}
\newcommand{\bed}[1]{\begin{eqnarray}\label{#1}}
\newcommand{\eed}{\end{eqnarray}}
\newcommand{\ind}[1]{\1_{\{ #1\}}}

\newcommand{\half}{\frac{1}{2}}

\newtheorem{theorem}{Theorem}[section]

\newcommand{\beds}{\begin{eqnarray*}}
\newcommand{\eeds}{\end{eqnarray*}}

\def\Z{{\mathcal Z}}

\def\V{\mathscr{V}} 
\def\A{\mathcal A}
\def\S{\mathbb S}
\def\T{\mathbb T}
\def\Normal{\mathscr{N}}

\def\Z{\mathscr{Z}}

\def\prob{\mathscr{P}_n}
\newcommand{\sel}[1]{\sup\Big\{ \prod_{j=1}^n n\pi_j: \bpi \in \prob,\  #1 \Big\}}
\newcommand{\isel}[1]{\sup\Big\{ \prod_{i=1}^n n\pi_i: \bpi \in \prob,\  #1 \Big\}}

\newcommand{\selN}[1]{\sup\Big\{ \prod_{j=1}^{N_n} N_n\pi_j: \bpi \in \mathscr{P}_{N_n},\  #1 \Big\}}

\def\H{\mathbb H}

\def\ep{{\epsilon}}
\def\vep{{\varepsilon}}

\def\h{\mathbf{h}}

\def\bdm#1{\mbox{\boldmath$#1$\unboldmath}}

\def\btheta{\mbox{\boldmath$\theta$\unboldmath}}

\def\bdelta{\bdm{\delta}}
\def\bpi{\bdm{\pi}}

\def\bvt{\bdm{\vartheta}}

\startlocaldefs
\numberwithin{equation}{section}

\numberwithin{equation}{section}
\theoremstyle{plain}
\def\Cov{\mathop{\rm Cov}\nolimits}
\def\Var{\mathop{\rm Var}\nolimits}%
\def\R{\mathbb R}
\def\x{\times}

\def\el{\mathscr{R}_{n}}

\def\prob{\mathscr{P}_{n}}
\def\und{\quad \mbox{and} \quad}

\def\1{\mathbf 1}
\def\Ave{\frac{1}{n}\sum_{j=1}^n}
\def\ave{n^{-1}\sum_{j=1}^n}

\def\iAve{\frac{1}{n}\sum_{i=1}^n}

\def\pave{\sum_{j=1}^n \pi_j}
\def\maxj{\max_{1\leq j \leq n}}

\def\trace{\mathrm{trace}}

\def\Z{{\mathbf Z}}

\def\V{\mathscr{V}}
\def\A{\mathcal A}
\def\S{\mathbb S}
\def\T{\mathbb T}

\def\U{\mathbb{U}_n}

\newcommand{\ssel}[1]{\sup\Big\{ \prod_{j=1}^n n\pi_j: &
(\pi_1,\dots,\pi_n)\in \Pi_n, \\ & #1 \Big\}}

\def\z{\mathfrak{z}}

\def\b{\mathbf{b}}
\def\u{\mathbf{u}}
\def\jsum{\sum_{j=1}^n}
\def\htet{\hbtheta}
\def\ttet{\tbtheta}
 
\def\seq#1{\left\{#1\right\}}
 
\def\jave{\frac{1}{n}\sum_{j=1}^n}

\def\lam{\lambda}

\def\seq#1{\left\{#1\right\}}
\def\V{\mathscr {V}}
\def\PS{Peng and Schick}
\def\jave{\frac{1}{n}\sum_{j=1}^n}

\def\pare#1{\left(#1\right)}

\def\Normal{\mathscr{N}}

\def\sZ{{\mathscr{Z}}}

\def\A{{\mathbf{A}}}
\def\B{{\mathbf{B}}}

\def\D{{\mathbf{D}}}

\def\H{{\mathbf{H}}}
\def\I{{\mathbf{I}}}
\def\J{{\mathbf{J}}}
\def\K{{\mathbf{K}}}

\def\M{{\mathbf{M}}}

\def\R{{\mathbf{R}}}
\def\T{{\mathbf{T}}}
\def\U{{\mathbf{U}}}
\def\V{{\mathbf{V}}}
\def\W{{\mathbf{W}}}
\def\X{{\mathbf{X}}}
\def\Y{{\mathbf{Y}}}
\def\Z{{\mathbf{Z}}}

\def\cR{{\mathcal{R}}}
\def\cS{{\mathcal{S}}}

\def\cZ{{\mathcal{Z}}}

\def\a{{\mathbf{a}}}
\def\b{{\mathbf{b}}}
\def\c{{\mathbf{c}}}

\def\f{{\mathbf{f}}}
\def\g{{\mathbf{g}}}
\def\h{{\mathbf{h}}}

\def\k{{\mathbf{k}}}

\def\r{{\mathbf{r}}}
\def\s{{\mathbf{s}}}
\def\t{{\mathbf{t}}}
\def\u{{\mathbf{u}}}
\def\v{{\mathbf{v}}}
\def\w{{\mathbf{w}}}
\def\x{{\mathbf{x}}}
\def\y{{\mathbf{y}}}
\def\z{{\mathbf{z}}}

\def\vs{{\otimes 2}}

\def\seq#1{\left\{#1\right\}}
\def\pare#1{\left(#1\right)}

\def\Cov{{\mathrm{Cov}}}

\def\jave{{\frac{1}{n}\sum_{j=1}^n}}

\def\bdm#1{\mbox{\boldmath$#1$\unboldmath}}
\def\btheta{\mbox{\boldmath$\theta$\unboldmath}}
\def\hbtheta{\bdm{\hat\btheta}}
\def\tbtheta{\bdm{\tilde\btheta}}

\def\bdelta{\bdm{\delta}}
\def\bpi{\bdm{\pi}}
\def\bzeta{\bdm{\zeta}}
\def\hbzeta{\bdm{\hat\zeta}}
\def\tbzeta{\bdm{\tilde\zeta}}

\def\bep{\bdm{\epsilon}}

\def\beeta{\bdm{\eta}}

\def\bvt{\bdm{\vartheta}}
\def\tbvt{\bdm{\tilde\vartheta}}
\def\bpsi{\bdm{\psi}}

\def\bphi{\bdm{\phi}}
\def\bvarphi{\bdm{\varphi}}
\def\bsigma{\bdm{\sigma}}
\def\bnu{\bdm{\nu}}
\def\hbvt{\bdm{\hat\vartheta}}
\def\bchi{\bdm{\chi}}
\def\bxi{\bdm{\xi}}

\def\vecs{\mathrm{vecs}}

\def\bmu{\bdm{\mu}}

\def\qed{$\square$}
\def\M{{\mathbf{\omega}}}
\def\tpi{{\tilde\pi}}

\def\bmS{\tilde{\bdm{\S}}}
\def\bms{\tilde{\bdm{\s}}}
\def\bPhi{\bdm{\Phi}}

\def\M{\mathbf{M}}
\def\qed{\hfill $\Box$}

\begin{document}
\renewcommand{\Ref}[1]{(\ref{#1})}

\begin{frontmatter}
\title{Improving Estimation Efficiency In Structural Equation Models 
By An Easy Empirical Likelihood Approach}
\runtitle{Improving Efficiency in SEM}

\begin{aug}
\author{\fnms{Shan}\snm{Wang} \thanksref{th1}
	\ead[label=e1]{swang151@usfca.edu}}
\and  
\author{\fnms{ Hanxiang} \snm{Peng} 
	\ead[label=e3]{hanxpeng@iu.edu}}
\thankstext{th1}{Corresponding author}
\runauthor{S. Wang and H. Peng}
	



\address{ University of San Francisco \\ 
         Department of Mathematics and Statistics\\
         San Francisco, CA 94117, USA \\
         swang151@usfca.edu
}

\address{ Indiana University Purdue University Indianapolis \\
         Department of Mathematical Sciences\\
         Indianapolis, IN 46202-3267, USA \\
         hanxpeng@iu.edu
}


\end{aug}

\begin{abstract}
In this article, we construct empirical likelihood (EL)-weighted estimators of linear functionals of a probability measure 
in the presence of side information. Motivated by nuisance parameters 
in semiparametric models with possibly infinite dimension, we consider the use of estimated constraint 
functions and allow the number of constraints to grow with the sample size. 
We study the asymptotic properties and efficiency gains.
The results are used to construct improved estimators of parameters in structural
equations models. 
The EL-weighted estimators of parameters are shown to have reduced variances 
in a SEM in the presence of side information of stochastic independence of the random error and random covariate. Some simulation results on efficiency gain are reported. 
\end{abstract}



\end{frontmatter}

\section{Introduction} \label{chsem}

Structure equation models (SEM) is  a popular multivariate technique for analyzing data in behavioral, medical and social sciences. It is an analysis of moment structures in which the variance-covariance matrix $\Sigma=\Var(\Z)$ of a random vector
$\Z\in \cR^p$ is specified by a parametric matrix function, $\Sigma=\Sigma(\bvt),\, 
\bvt \in \Theta$ for some subset $\Theta$ of $\cR^q$.  Given independent and
identically distributed (i.i.d.) observations $\Z_1, \dots, \Z_n$ of $\Z$, one focuses on estimating 
the parameter vector $\bvt$.  A moment-type estimator of $\bvt$ is based on the criterion of 
\emph{minimum discrepancy function (MDF)}, see e.g. Shapiro (2007)\cite{sh07}. 
For a $p \times p$ matrix $\M$, denote by $\vecs(\M)$ the $p(p+1)/2$-dimensional vector formed by stacking its columns
of the upper triangular matrix.
Let $\Xi$  be a subset of $\cR^{p(p+1)/2}$ consisting of $\vecs(\M)$  over all  $p \times p$ semi-positive definite
matrices $\M$. A function $F$ on $\Xi \times \Xi$ is called a \emph{discrepancy function} if it satisfies:  
(i) $F(\t, \bxi)\geq 0$ for all $\t, \bxi \in \Xi$,
(ii) $F(\t, \bxi)=0$ if and only if $\t=\bxi$,
(iii) $F(\t, \bxi)$ is twice continuously differentiable, and  
(iv) For any fixed $\s \in \Xi$, $F(\t, \bxi) \to \infty$  as $\t \to \s$ and 
	$\bxi\to\infty$. 
We shall abuse notation to write $F(\A, \B)=F(\vecs({\A}), \vecs(\B))$ for matrices $\A, \B$.

  Often $\Sigma$ is estimated by the unstructured sample variance-covariance matrix,   
\bel{mst1}
\S_n=\Avi (\Z_i-\bar \Z)(\Z_i-\bar \Z)^\top,
\ee
where $\bar \Z=\Avi \Z_i$ is the sample version of $\bmu=E(\Z)$. Let $\s_n=\vecs(\S_n)$ and $\bsigma=\vecs(\Sigma)$. The MDF estimator  of $\bvt$ is any value $\bvt_n$ in $\Theta$ that satisfies
\bel{mst12}
F(\s_n, \sigma(\bvt_n))=\inf_{\bvt \in \Theta} F(\s_n, \bsigma(\bvt)).
\ee
 One commonly used discrepancy function is the maximum-likelihood (ML) discrepancy function of two matrix-valued variables given by
\bel{ml}
F_{ML}({\mathbf S}, \Sigma)=\log |\Sigma| -\log |{\mathbf S}|
+\trace({\mathbf  S} \Sigma^{-1})-p, 
\ee
where $\mathbf S, \Sigma$ are positive definite matrices. 
Another is the generalized least squares (GLS) discrepancy
function given by 
\bel{gls}
F_{GLS}(\mathbf {S}, \Sigma)=\trace((\mathbf {S}-\Sigma){\W}^{-1}({\mathbf S}-\Sigma){\W}^{-1}),  
\ee
where ${\W}$ is a symmetric matrix, for instance, ${\W}=\S_n$ and ${\W}=\I_p$. 
 Suppose that there is available some information about SEM that can be expressed by a vector equation of expectation,  
\bel{sinfo}
E(\g(\Z))=0, 
\ee
where $\g$ is  some measurable function taking values in $\cR^m$. 
While the moment estimator $\S_n$ completely ignores  the information, the empirical likelihood (EL)-weighted estimator $\bmS_n$  utilizes the information to provide an improved estimator, 
\bel{mst2}
\bmS_n=\Avi \frac{(\Z_i-\bar \Z)(\Z_i-\bar \Z)^\top}{1+\g(\Z_i)^\top \tbzeta}, 
\ee
where $\tbzeta$ is a solution to the equation
\bel{mst21}
\sum_{i=1}^n \frac{\g(\Z_i)}{1+\g(\Z_i)^\top \tbzeta}=0. 
\ee
Accordingly, an improved MDF estimator $\tbvt_n$ of $\bvt$ is any value in $\Theta$ that satisfies
\bel{mst22}
F(\bms_n, \sigma(\tbvt_n))=\inf_{\bvt \in \Theta} F(\bms_n, \sigma(\bvt)),
\ee
where $\bms_n=\vecs(\bmS_n)$. 

In many semiparamatric models, $\g(\z)$ involves in nuisance parameters which must be estimated, 
leading to a plug-in estimator $\hat\g(\z)$. Using it, we work with 
\bel{mst23}
\hat\S_n=\Avi \frac{(\Z_i-\bar \Z)(\Z_i-\bar \Z)^\top}{1+\hat\g^\top(\Z_i) \hbzeta}, 
\ee
where $\hbzeta$ is the solution to Eqt \Ref{mst21} by replacing $\g(\Z_i)=\hat\g(\Z_i)$.
As a result,  an improved  MDF estimator $\hbvt_n$ of $\bvt$ is any value in $\Theta$ that satisfies
\bel{mst25}
F(\hat\s_n, \sigma(\hbvt_n))=\inf_{\bvt \in \Theta} F(\hat\s_n, \sigma(\bvt)), 
\ee
where $\hat\s_n=\vecs(\hat\S_n)$.

  The improved MDF estimator $\tbvt_n$ is more efficient than the usual MDF estimator $\bvt_n$. 
The efficiency criteria used are that of a least dispersed regular estimator or that of a locally asymptotic minimax estimator, and are based on the convolution theorems and
on the lower bounds of the local asymptotic risk in LAN and LAMN families,
see the monograph by Bickel, {\sl et al.} (1993)\cite{bk93}. 

The side information contained in \Ref{sinfo} is carried by the EL-weights $(n(1+\g(\Z_i)^\top \tbzeta))^{-1}$ based on the principle of the maximum empirical likelihood.  
There is an extensive amount of literature on the empirical likelihood. 
It was introduced by Owen (1990, 1991) \cite{O90, O91} to construct  confidence intervals in a nonparametric setting.
Soon it was used to construct point estimators. Qin and Lawless (1994) \cite{QL94} studied 
maximum empirical likelihood estimators (MELE). 
Bravo (2010)\cite{br10} studied a class of M-estimators based on generalized empirical likelihood  
 with side information and showed that the resulting class of 
estimators is efficient in the sense that it achieves the same asymptotic lower bound as that of 
the efficient GMM estimator with the same side information.  
   Parente and Smith (2011)\cite{PS11}  investigated generalized empirical likelihood estimators 
for irregular constraints.

  Estimators of the preceding EL-weighted form were investigated in Zhang (1995, 1997) \cite{z95, z97} 
in M-estimation and quantile processes in the presence of auxiliary information.   
Hellerstein and Imbens (1999) \cite{hi99} exploited such estimators  for the least squares estimators 
in a linear regression model.   
Yuan {\sl et al.} (2012) \cite{yh12} explored  such estimators in U-statistics. 
Tang and Leng (2012) \cite{tl12} utilized the form to construct improved estimators of 
parameters in quantile regression. Asymptotic properties of the EL-weighted estimators 
were obtained for a {finite} number of {known} constraints.  Motivated by nuisance parameters in semiparametric models and the infinite dimension of such models, 
Peng and Schick (2013)\cite{ps13} considered the use of estimated constraint functions and studied  
 a growing number of constraints in MELE. MELE enjoy high efficiency and is particularly convenient to incorporate side information.
Just like any other optimization problems, however, it is not trivial to numerically find MELE especially 
for a large number of constraints. 
Peng and Schick \cite{ps13} employed one-step estimators to construct MELE.
%
The EL-weighted approach reduces the number of constraints 
and are thus computationally easier than general MELE. 

   This article used the EL-weighted approach to construct 
efficient estimators of linear functionals of a probability measure in the presence of side information 
for two cases, viz,  known marginal distributions and equal but unknown marginals, 
each of which is equivalent to infinitely many constraints. The rest of the article is organized as follows. 
In Section \ref{main}, we shall construct the EL-weighted estimator of the linear functional of a probability
measure in the presence of side information which is expressed by an finite or infinite number of known 
or estimated constraints, and present the asymptotic properties. 
 In Section \ref{semmainthm}, we give examples of side information and study 
the asymptotic properties of the improved estimators in SEM. 
 The form of SEM can be extended to great extent in a variety of ways.  
We shall focus on the extensions that have been described in Bollen (1989)\cite{bo89} as well as in the LISREL software manual (Jöreskog and Sörbom (1996)\cite{js96}). 
The components present in a general SEM are a path analysis, the conceptual synthesis of latent variable and measurement models, and general estimation procedures. In SEM, only information up to the second moments 
is used, while other forms of information such higher order moments, independence or symmetry 
of the random errors are ignored, which can be used by the EL-weighting method to improve efficiency.
In Section \ref{sim}, we report simulation results. Technical details are collected in Section \ref{proofs}.

\section{The main results}\label{main} 
Suppose that $Z_1, \dots, Z_n$ are i.i.d. random variables with a common distribution $Q$ taking values in a measurable space $\cZ$. 
We are interested in efficient estimation of the linear functional 
$\btheta=\int \bpsi\,dQ$ of $Q$ for some square-integrable function $\bpsi$ from $\cZ$ to $\cR^r$ when side information is available through 
\begin{itemize}
	\item[(C)] $\u$ is a measurable function from $\cZ$ to $\cR^m$ such that $\int \u\,dQ=0$ and
	the variance-covariance matrix $\W=\int\u\u^\top\,dQ$ is nonsingular.
\end{itemize}
To utilize the information contained (C), consider the the empirical likelihood, 
$$
\el=\sel{ \pave \u(Z_j)=0},  
$$
where $\prob=\{\pi \in [0, \, 1]^n: \sum_{j=1}^n \pi_j=1\}$ is the unit probability simplex. 
Following Owen (1990)\cite{O90}, one uses Lagrange multipliers to get  the maximizers, 
\bel{elw}
\tpi_{j}=\frac{1}{n}\frac{1}{1+\u(Z_j)^\top\tbzeta}, \quad j=1, \dots, n,
\end{equation}
where $\tbzeta$ is the solution to the equation  
\bel{tele2}
\Ave \frac{\u(Z_j)}{1+\u(Z_j)^\top\tbzeta}=0.
\end{equation}  
These $\tpi_{j}$'s incorporate the side information, and a natural estimator of $\btheta=\int\bpsi\,dQ$  
is the EL-weighted estimator, 
\bel{tele}
\tbtheta=\sum_{j=1}^n \tpi_{j}\bpsi(Z_j)=\Ave \frac{\bpsi(Z_j)}{1+\u(Z_j)^\top \tbzeta}.  
\ee 
For $\bpsi_\t(\z)=\1[\z\leq \t]$ for fixed $\t \in\cR^p$, one obtains the distribution function
$\theta=P(\Z\leq \t)$. For $\psi(\z)=z_1\cdots z_p$, $\theta=E(Z_1\cdots Z_p)$ is the mixed
moment.

Write $\|\a\|$ for the euclidean norm of $\a$ and $\a \otimes \b$ for the Kronecker product of $\a$ and $ \b$.
For $\x=(x_1, \dots, x_p), \y=(y_1, \dots, y_p)$, write $\x\leq \y$ for $x_1\leq y_1, \dots, x_p\leq x_p$.  
Let $L_2^m(Q)=\seq{\f=(f_1, \dots, f_m)^\top: \int \|\f\|^2\,dQ^m<\infty}$, and let
$L_{2,0}^m(Q)=\seq{\f\in L_2^m(Q): \int \f\,dQ^m=0}$. 
For $\f \in L_2^m(Q)$, write $[\f]$ for the closed linear span of the components $f_1$, $\dots$, $f_m$ in $L_2(Q)$. 
Let $Z$ be an i.i.d. copy of $Z_1$. 
Let $\bphi_0$ be the projection of $\bpsi$ onto the closed linear span $[\u]$
of $\u$, so that $\bphi_0=\Pi(\bpsi|[\u])=E(\bpsi(Z)\otimes\u^\top(Z))\W^{-1}\u$. 
Let $\Sigma_0=\Var(\bpsi(Z))-\Var(\bphi_0(Z))$. We now give the first result with the proof delayed in Section \ref{proofs}.  

\Thm{1} 
Assume (C) with $m$ fixed. 
Then $\ttet$ given in \Ref{tele} satisfies the stochastic expansion, 
\bel{al1}
\ttet=\bar \bpsi- \bar \bphi_0+o_p(n^{-1/2}),
\end{equation} 
Thus if $\Sigma_0=\Var(\bpsi(Z))-\Var(\bphi_0(Z))$ is nonsingular then 
$\sqrt{n}(\ttet-\btheta)$ is asymptotically normal with mean zero and asymptotic covariance matrix $\Sigma_0$, that is, 
$$
\sqrt{n}(\ttet-\btheta) \implies \Normal(0, \Sigma_0).
$$

\end{theorem} 

\thmref{1} exhibits that the EL-weighted estimator $\ttet$ has a smaller asymptotic variance than
that of the sample mean $\bar\bpsi$, and the amount of reduction is $\Var(\bphi_0(Z))$.  It  is, in fact, the MELE of $\btheta$.

\Rem{hab} Haberman (1984)\cite{h84} studied minimum Kullback-Leibler divergence -type estimators for the linear functionals
of a probability measure, and more general problems involving a fixed number of side information. 
The EL-weighted estimator $\ttet$ in \thmref{1} is asymptotically equivalent to Haberman's estimator, 
see his page 976. 
This shows that Haberman's estimator is semiparametrically efficient. 
\end{remark}

  In semiparametric models, the constraint function $\u$ contains nuisance parameters and must be estimated. 
Let $\hat\u=(\hat u_1, \dots, \hat u_m)^\top$ be an estimate of $\u$. With it we now work with the EL-weights  
\bel{ho1}
\hat\pi_j=\frac{1}{n}\frac{1}{1+\hat\u(Z_j)^\top \hbzeta}, \quad j=1,\dots,n,
\ee
where $\hbzeta$ is the solution to the equation \Ref{tele2} with $\u=\hat\u$.   
In the same fashion, a natural estimate $\htet$ of $\btheta$ is given by 
\bel{hele}
\htet=\jsum \hat\pi_j\bpsi(Z_j)=\Ave \frac{\bpsi(Z_j)}{1+\hat \u(Z_j)^\top \hbzeta}.
\end{equation}   
Set $\hat \W=\ave \hat\u\hat\u^\top(Z_j)$. Let $|\W|_o$ denote the spectral norm (largest eigenvalue) of a matrix $\W$.
We have  

\Thm{2} 
Assume (C) with $m$ fixed.  Let $\hat \u$ be an estimator of $\u$ such that
\bel{2a}
\maxj \|\hat\u(Z_j)\|=o_p(n^{1/2}), 
\end{equation}
\bel{2b}
|\hat \W-\W|_o=o_p(1), 
\end{equation}
\bel{2c}
\Ave \pare{\bpsi(Z_j)\otimes \hat\u(Z_j)-E\big(\bpsi(Z_j)\otimes \hat\u(Z_j)\big)}=o_p(1),
\end{equation}
and  that there exists some measurable function $\v$ that satisfies (C) such that    
\bel{2d}
\Ave E\pare{\|\hat\u(Z_j)-\v(Z_j)\|^2}=o(1),
\end{equation}
\bel{2e}
\Ave \hat \u(Z_j)=\Ave\v(Z_j)+ o_p(n^{-1/2}). 
\ee
Then $\htet$ given in \Ref{hele} satisfies  the stochastic expansion, 
\bel{al2}
\htet=\bar \bpsi-\bar\bphi+o_p(n^{-1/2}),
\end{equation} 
where $\bphi=\Pi(\bpsi|[\v])$. 
Thus if $\Sigma=\Var(\bpsi(Z))-\Var(\bphi(Z))$ is nonsingular then
$$
\sqrt{n}(\htet-\btheta) \implies \Normal(0, \Sigma). 
$$

\end{theorem} 
\def\Nu{{\mathcal N}}
\def\hbvt{\bdm{\check\btheta}}

\def\hbvt{\bdm{\hat\bvt}}

We now allow the number of constraints to depend on $n$, 
$m=m_n$,  and grow to infinity with increasing $n$. To stress the dependence, let us
write 
$$\u_n=(u_1, \dots, u_{m_n})^\top, \quad
\hat \u_n=(\hat u_1, \dots, \hat u_{m_n})^\top, 
$$ 
and $\ttet_n=\ttet$, $\htet_n=\htet$ for 
the corresponding estimators of $\btheta$, that is,
\bel{telen}
\ttet_n=\Ave \frac{\bpsi(Z_j)}{1+\u_n(Z_j)^\top \tilde\bzeta_n} \und
\htet_n=\Ave \frac{\bpsi(Z_j)}{1+\hat\u_n(Z_j)^\top \hat\bzeta_n},
\end{equation}
where $\tilde\bzeta_n$ and $\hat\bzeta_n$ solves Eqt 
\Ref{tele2} with $\u=\u_n$ and $\u=\hat\u_n$, respectively. 
Denote by $[\u_\infty]$ the closed
linear span of $\u_\infty=(u_1, u_2, \dots)$. Set 
$$ 
\W_n=\Var(\u_n(Z)), \quad
\bar \W_n=\Ave \u_n\u_n^\top(Z_j), \quad
\hat \W_n=\Ave \hat\u_n\hat\u_n^\top(Z_j).
$$

 Peng and Schick (2013) \cite{ps13} introduced that a sequence $\W_n$ of $m_n\times m_n$ dispersion matrices 
is \emph{regular} if
$$
0 < \inf_n \inf_{\|\u\|=1} \u^{\top}\W_n \u \leq \sup_n \sup_{\|\u\|=1} \u^{\top}\W_n \u <\infty.
$$
  Note that if $\W=\W_n$ is independent of $n$ then the regularity of $\W$ simplifies to its nonsingularity. 
We have 

\Thm{3} 
Suppose that $\u_{n}=(u_1, \dots, u_{m_{n}})^\top$ satisfies (C) for each $m=m_n$  
such that
\bel{3a}
\maxj \|\u_n(Z_j)\|=o_p(m_n^{-3/2}n^{1/2}), 
\ee 
that the sequence of $m_n \times m_n$ dispersion matrices $\W_n$ is {regular} and satisfies
\bel{3b}
|\bar \W_n-\W_n|_o=o_p(m_n^{-1}), 
\ee
\bel{3c}
\Ave\pare{\bpsi(Z_j)\otimes \u_n(Z_j)-E\big(\bpsi(Z_j)\otimes \u_n(Z_j)\big)}=o_p(m_n^{-1/2}).
\ee
Then
$\ttet_n$ satisfies, as $m_n$ grows to infinity with $n$,  the stochastic expansion, 
\bel{al3}
\ttet_n=\bar \bpsi-\bar\bvarphi_0+o_p(n^{-1/2}),
\end{equation} 
where $\bvarphi_0=\Pi(\bpsi|[\u_\infty])$.  
Thus if $\varSigma_0=\Var(\bpsi(Z))-\Var(\bvarphi_0(Z))$ is nonsingular, 
$$
\sqrt{n}(\ttet_n-\btheta) \implies \Normal(0, \varSigma_0). 
$$

\end{theorem} 

\Thm{4} Suppose that $\u_{n}=(u_1, \dots, u_{m_{n}})^\top$ satisfies (C) for each $m=m_n$.    
Let $\hat\u_n$ be an estimator of $\u_n$  such that
\bel{4a}
\maxj \|\hat\u_n(Z_j)\|=o_p(m_n^{-3/2}n^{1/2}), 
\end{equation}
\bel{4b}
|\hat \W_n-\W_n|_o=o_p(m_n^{-1})
\end{equation}
for which the $m_n \times m_n$  dispersion matrices $\W_n$ is regular, 
\bel{4c}
\Ave\pare{\bpsi(Z_j)\otimes \hat\u_n(Z_j)-E\big(\bpsi(Z_j)\otimes \hat\u_n(Z_j)\big)}=o_p(m_n^{-1/2}),
\end{equation}
and that there exists some measurable function $\v_n$ from $\cZ$ into $\cR^{m_n}$ such that 
(C) is met for every $m=m_n$, the dispersion matrix $\U_n=\W_n^{-1/2}\int \v_n\v_n^\top\,dQ\W_n^{-\top/2}$ satisfies $|\U_n|_o=O(1)$, and 
\bel{4d}
\Ave E\pare{\|\hat\u_n(Z_j)-\v_n(Z_j)\|^2}=o(m_n^{-1}), \quad \textrm{and} 
\end{equation}
\bel{4e}
\Ave \hat \u_n(Z_j)=\Ave\v_n(Z_j)+ o_p(m_n^{-1/2} n^{-1/2}). 
\ee
Then $\htet$ satisfies, as $m_n$ tends to infinity,  the stochastic expansion, 
\bel{al4}
\htet_n=\bar \bpsi-\bar\bvarphi+o_p(n^{-1/2}),
\end{equation} 
where $\bvarphi=\Pi(\bpsi|[\v_\infty])$. 
Thus if $\varSigma=\Var(\bpsi(Z))-\Var(\bvarphi(Z))$ is nonsingular,  
$$
\sqrt{n}(\htet_n-\btheta) \implies \Normal(0, \varSigma).
$$
\end{theorem}

\section{The asymptotic properties of EL-weighted MDF estimators}\label{semmainthm}
In this section, we give two examples of side information and discuss  
the reduction in the asymptotic covariance matrix 
of the improved MDF estimator $\hat\bvt_n$  given 
in \Ref{mst25} compared with the sample MDF estimator $\bvt_n$.

\subsection{Examples and side information} We introduce the SEM and discuss side information.
\Ex{sem} Consider the \emph{combined model} of latent variable and measurement error, 
\bel{sem}
\beeta=\
\B\beeta+\Gamma \bxi +\bzeta, \quad \Y-\bmu_y=\Lambda_y \beeta+\bep, \quad \X-\bmu_x=\Lambda_x\bxi+\bdelta, 
\ee
where 
$\B$, $\Gamma$, $\Lambda_x$, $\Lambda_y$, $\bmu_x$ and $\bmu_y$ are compatible parameter matrices
and vectors, $\X$ and $\Y$ are random vectors having finite fourth moments, 
$\beeta$ and $\bxi$ are  latent endogenous and exogenous random vectors,  respectively, 
and $\bzeta$, $\bep$ and $\bdelta$ are disturbances (random errors) that satisfy 
%
\bel{sem1}
\begin{aligned}
	& E(\bzeta)=0,\; E(\bep)=0, \; E(\bdelta)=0, 
	\quad \Cov(\bep, \beeta)=0,
	\quad \Cov(\bdelta, \bxi)=0,\\
	& \Cov(\bxi, \, \bzeta)=0,\quad 
	\Cov(\bep, \, \bzeta)=0,\quad 
	\Cov(\bdelta, \, \bzeta)=0, \quad 
	\Cov(\bep, \, \bdelta)=0.  
\end{aligned} 
\ee
Let $\Phi=E(\bxi\bxi^\top)$, $\Psi=E(\bzeta\bzeta^\top)$,  $\Theta_\ep=E(\bep\bep^\top)$
and $\Theta_\delta=E(\bdelta\bdelta^\top)$. The parameter vector then is $\bvt=\vecs(\bmu_x, \bmu_y, \B, \Gamma, \Lambda_x, \Lambda_y, \Phi, \Psi, \Theta_\ep, \Theta_\delta)$, denoted by $q$ the dimension. 
Let $\Sigma_{yy}(\bvt)$ be the structured variance-covariance of $\Y$, and let $\A=\I_d-\B$. Based on the relationships in \Ref{sem} -- \Ref{sem1}, one derives  
$$
\Sigma_{yy}(\bvt)=\Lambda_y\A^{-1}(\Gamma\Phi\Gamma^\top+\Psi)\A^{-\top}\Lambda_y^\top+\Theta_\ep,
$$
 assuming that $\A$ is invertible. Similarly, one derives the structured covariance matrix $\Sigma_{yx}(\bvt)$ of $\Y$ and $\X$ and the variance-covariance $\Sigma_{xx}(\bvt)$ of $\X$, 
$$
\Sigma_{yx}(\bvt)=\Lambda_y\A^{-1}\Gamma\Phi\Lambda_x^\top=\Sigma_{xy}(\bvt)^\top, \quad
\Sigma_{xx}(\bvt)=\Lambda_x\Phi\Lambda_x^\top+\Theta_\delta.
$$
The structured variance-covariance $\Sigma(\bvt)$ of $\Z=(\Y^\top, \X^\top)^\top$ then is
$$
\Sigma(\bvt)=\begin{pmatrix}
\Sigma_{yy}(\bvt) & \Sigma_{yx}(\bvt)\\
\Sigma_{xy}(\bvt) & \Sigma_{xx}(\bvt)
\end{pmatrix}. 
$$
These formulas can be found in literature, but we would mention that they are implied 
by the structural relationships in \Ref{sem} -- \Ref{sem1}. 
While the unstructured sample variance-covariance matrix estimator $\S_n$ in \Ref{mst1}
of the unstructured variance-covariance $\Sigma$ of  $\Z$ 
ignores the information contained in \Ref{sem1},  the EL-weighted 
estimator $\hat\S_n$ of $\Sigma$ in \Ref{mst23} utilizes the information, 
and results in an improved estimator $\hbvt$ determined by \Ref{mst25}. 
\end{example}

\Ex {sems}  In the combined model in \exref{sem}, consider $\Lambda_y=\I_d$, $\Lambda_x=\I_c$,
$\Var(\bdelta)=0$ and $\Var(\bep)=0$. This is a SEM, and \Ref{sem} -- \Ref{sem1} simplify to 
\bel{sems1}
\Y=\B\Y+\Gamma \X +\bzeta,  \quad E(\bzeta)=0, \quad \Cov(\X, \, \bzeta)=0.
\ee
    Identification is crucial for the consistency and asymptotic normality of the MDF estimators. Necessary and sufficient conditions can be found in the literature, e.g.,  Bollen (1989)\cite{bo89}, 
Brito and Pearl (2002)\cite{bp02} and Drton, {\sl et al.} (2011)\cite{dr11}.
In particular, the Null B Rule and the Recursive Rules are sufficient conditions for the identifiability
of the parameters. The former states that if $\B=0$ then the parameters can be identified, while the
latter says that if $\B$ can be written as a lower triangular matrix with zero diagonal and the covariance matrix $\Psi$ of the error $\bzeta$ is diagonal then the parameters are identifiable.  
An example for the latter case is the model given by 
\bel{sems2}
\Y=\B\Y+\Lambda \X+\bep, 
\ee
where $\B, \Lambda$ are $2\times 2$ matrices, with $\B$ having all entries equal to 0 except for the $(2, 1)$ entry equal to $\beta$, and $\Lambda$ having the (1, 1) entry  equal to 0 and the (1, 2), (2, 1) and (2, 2) entries equal to $\lambda_i, i=1, 2, 3$, respectively. The path diagram is shown in Fig. 1.
\begin{figure} \label{figure_1}
	\begin{center}
		\newcommand\nodelabel[3][3em]{\parbox{#1}{\centering #2 \\ #3 }}
		\newcommand\vertspacing{-1.5}
		\begin{tikzpicture}[align=center,node distance=3cm,
		]
		\tikzstyle{paths}=[-latex,black,line width=1pt]
		\tikzstyle{var}=[latex-latex,black,line width=1pt]
		\tikzstyle{cov}=[latex-latex,black,line width=1pt,dashed, bend right=270, looseness=0.8]
		\tikzstyle{obs}=[rectangle,draw,font=\scriptsize\sffamily\bfseries,inner sep=3pt]
		\tikzstyle{latent}=[ellipse,draw, font=\scriptsize\sffamily\bfseries,inner sep=3pt]
		\tikzstyle{error}=[circle,draw,
		font=\scriptsize\sffamily\bfseries,inner sep=3pt]
		\node [error](e2) at (0,0) { $\varepsilon_2$};
		\node [error](e1) [left of=e2] { $\varepsilon_1$};
		\node [obs](y1) [below =1cm of e1]{ $y_1$};
		\node [obs](y2) [below=1cm of e2]{$y_2$};
		\node [obs](x1) [below left=1cm and 0.5cm of y1]{ $x_1$};
		\node [obs](x2) [below right=1cm and 0.8cm of y2]{ $x_2$};
		\draw [paths,above](e1) to (y1);
		\draw [paths,above](e2) to (y2);
		\draw [paths,above](y1) to node [above=0.3cm] {$\beta$} (y2);
		\draw [paths,above](x2) to node [below=0.1cm] {$\lambda_1$} (y1);
		\draw [paths,above](x1) to node [left=0.4cm] {$\lambda_2$} (y2);
		\draw [paths,above](x2) to node [right=0.25cm] {$\lambda_3$} (y2);
		\end{tikzpicture}
	\end{center}
	\caption{The path diagram for SEM \Ref{sems2} }
\end{figure}

{\bf Side information}. 
SEM make use of the information up to the second moments, whereas other information is completely ignored.  
For example, random errors are modeled as
uncorrelated with covariates. It is common that the random error $\ep$ is modeled as independent of the
random covariate $\X$. The information contained in the independence can be utilized by  
the vector constraint function,
\bel{sems3}
\g(\Z)=\bPhi_m(F(\vep))\otimes \bPhi_m (G(\X)),  
\ee
where $\bPhi_m(t)=\sqrt{2}(\cos(\pi t), \dots, \cos (m \pi t))^\top$ is a vector of the first $m$ terms of
the trigonometric basis, and $F$ and $G$ are the respective distribution functions (DF) of the linear combination 
$\vep=\a^\top\bep$ of $\bep$ and $\X$. Here $\a$ is a known constant vector and  
$\otimes$ denotes the Kronecker product.
See Example 1 of Peng and Schick (2013)\cite{ps13} for more details.  
As $F, G$ are unknown, we estimate them by the empirical DF (EDF)   
$F_n, G_n$. We replace $\bep$ with 
$\hat{\bep}=\Y-\hat\B \Y-\hat{\Gamma} \X$,  where $\hat\B$ and  
$\hat {\Gamma}$ are  the MDF estimators of $\B$ and $\Gamma$. 
Substitution of them in \Ref{sems3} yields the estimated constraint function, 
\bel{sems4}
\hat\g(\Z)=\bPhi_m(F_n(\hat\vep))\otimes \bPhi_m (G_n(\X)),  
\ee 
   This is a  semiparametric model with (infinite dimensional) nuisance parameters $F, G$, and  
the plug-in estimators of  $F_n, G_n$ lead to the estimated constraints.

 Another example of side information is that the marginal medians (or means)  
$m_{01}$ and $m_{02}$ of $\X$ are  \emph{known}. Such marginal information is often possible
such as from the past data.  
In this case,  
the constraint function is 
\bel{sems5}
{\g}(\Z_j)=(\1[X_{1j} \leq m_{01}]-0.5, \, \1[X_{2j}\leq m_{02}]-0.5)^\top, \,
j=1, \dots, n.
\ee      

\end{example}


\subsection{The asymptotic properties}
We need some results from Shapiro (2007)\cite{sh07}. 
Let $\bvt_0$ be the true value of parameter $\bvt$ and $\bxi_0=\bsigma(\bvt_0)$.
By the Taylor expansion it is not difficult to show that a discrepancy function $F$ satisfies
\bel{mst3}
2\H_0:=\frac{\partial^2 F(\bxi_0, \bxi_0)}{\partial\t\partial\t^\top}
=\frac{\partial^2 F(\bxi_0, \bxi_0)}{\partial\bxi\partial\bxi^\top}
=-\frac{\partial^2 F(\bxi_0, \bxi_0)}{\partial\t\partial\bxi^\top},
\ee
and $\H_0$ is positive definite, see also Shapiro (2007)\cite{sh07}. 
In particular, for both $F_{ML}$ and $F_{GLS}$ (in the case of $\sl W=\S_n$), one has 
\bel{mst4}
\H_0=\Sigma_0^{-1} \otimes \Sigma_0^{-1}, 
\ee
where $\Sigma_0=\Sigma(\bvt_0)$.
Formally, set $\Delta(\bvt)=\partial \bsigma(\bvt)/\partial \bvt^\top$  with $\Delta_0=\Delta(\bvt_0)$ and 
\bel{mst5}
\begin{aligned}
	& \w(\z)=\vecs\big((\z-\bmu_0)(\z-\bmu_0)^\top-\Sigma_0\big), \quad  \z\in\cR^p, \\
	& \v(\z)=\g(\z)+E(\dot \g(\Z))\Psi(\z), \quad \Psi(\z)=(\Delta_0^\top \H_0 \Delta_0)^{-1}\Delta_0^\top \H_0 \w(\z).   
\end{aligned}
\ee
Summarizing Shapiro's results, we have 

\Lm{mst} Let $\Z, \Z_1, \dots, \Z_n$ be i.i.d. random vectors with finite and nonsingular
covariance matrix  
$\Var(\w(\Z))$. Assume that $\bvt_0$ is an interior point of $\Theta$ which is compact and can be approximated at $\bvt_0$ by $\cR^q$.  Suppose that $F$ is a discrepancy function. Suppose that $\bsigma(\bvt)$ is twice continuously differentiable with gradient $\Delta(\bvt)$ of full rank $q$ in a neighborhood of $\bvt_0$. 
Suppose that the model is locally identifiable, i.e.,
$\bsigma(\bvt)=\bsigma(\bvt_0)$ implies $\bvt=\bvt_0$ for $\bvt$ in a neighborhood of $\bvt_0$.
Then 
\bel{mst63}
\sqrt{n}(\tbvt-\bvt_0)\implies \Normal(0, \V_0),
\ee
where $\V_0=(\Delta_0^\top \H_0 \Delta_0)^{-1}\Delta_0^\top \H_0\Var(\w(\Z))
\H_0 \Delta_0(\Delta_0^\top \H_0 \Delta_0)^{-1}$. 

\end{lemma}

\Rem{resi} \lmref{mst} implies $\tbvt-\bvt_0=O_p(n^{-1/2})$. 
Consequently, each residual satisfies $\hat\bep_i-\bep_i=O_p(n^{-1/2})$ for  $\Z_i$ of bounded second
moment. 
We shall impose a stronger assumption of $E\|\hat\bep_i-\bep_i\|^2=O(n^{-1})$ uniformly in $i$. 
\end{remark}

{\sc Proof of \lmref{mst}}. We shall present the proof based on Theorem 5.5 of Shapiro (2007)\cite{sh07}. 
To this end, we first verify the conditions of his Proposition 4.2 to show $\tbvt$ is a consistent estimator of $\bvt_0$. 
Note that $\s_n=\vecs(\S_n)$ is clearly a (strongly) consistent estimator of $\bsigma_0=\bsigma(\bvt_0)$ since
$\S_n$ is a (strongly) consistent estimator $\Sigma_0=\Sigma(\bvt_0)$.
The local identifiability of $\bsigma(\bvt)$ at $\bvt_0$ implies the uniqueness of the optimal
solution (i.e. $\bvt_0$), hence his (4.4) is proved since $\Theta$ is compact, see the last paragraph 
of his page 238. This establishes the consistency by his Proposition 4.2.
It thus follows from his Theorem 5.5, (5.11) and (5.33) that $\tbvt$ satisfies
\bel{mst632}
\tbvt=\bvt_0+(\Delta_0^\top \H_0 \Delta_0)^{-1}\Delta_0^\top \H_0(\t_n-\bsigma_0)+o_p(n^{-1/2}),
\ee
where $\t_n=\vecs(\T_n)$ with $\T_n=\ave(\Z_j-\bmu_0)(\Z_j-\bmu_0)^\top$. Since $\Z$ has finite fourth moment,  it follows from the central limit theorem,  
\bel{mst64}
\sqrt{n}(\t_n-\bsigma_0) \implies \Normal(0, \Var(\w(\Z))).
\ee
The preceding two displays yield the desired \Ref{mst63} and end  the proof. 
\qed

R\'ev\'esz (1976)\cite{rp76} investigated the approximation of the empirical distribution function in two dimension.
Unlike in the case of one dimension in which the Kolmogorov-Smirnov statistic is asymptotically distribution free, 
the test in two dimension is not asymptotically distribution free, as shown in his Theorem 3,
which is quoted in \lmref{r1} below.
 Let $\Y=(Y_{1}, Y_{2})^\top$ be a random vector, and let $\T$ be a transformation 
of $\Y$ on $\cR^2$ such that $\T\Y$ is uniformly distributed. 
Consider the transformation given by  
$\T(y_1, y_2)=(H(y_1), G(y_2|y_1))^\top$, where  
\bel{condt}
H(y_1)=P(Y_{1} \leq y_1), \quad 
G(y_2|y_1)=P(Y_{2} \leq y_2|Y_{1}=y_1).
\ee

\Lm{r1}
Let $\Y_1=(Y_{11},Y_{12})$, $\Y_2=(Y_{21},Y_{22})$, $\dots$ be a sequence of i.i.d. rv's having a common DF  $F(\y)=F(y_1, y_2)$. Suppose that $F(y_1, y_2)$ is absolutely continuous and satisfies  
\bel{re1}
\Big|\frac{\partial G(y_2|H^{-1}(y_1))}{\partial y_1}\Big| \leq L, \,
\Big| \frac{\partial^2 G(y_2|H^{-1}(y_1))}{\partial y_1^2}\Big|\leq L,  
\, \y=(y_1, y_2) \in \cR^2, 
\ee
for some constant $L>0$.
Then we can define a sequence $\{\bar{B}_n\}$ of Brownian Measures (B.M.) and a Kiefer Measure (K.M.) 
$\bar{K}$ such that 
\bel{mks}
\begin{aligned}
&\sup_{\y\in \cR^2} |\beta_n(\y)-\bar{B}_n(TD_\y)|=O(n^{-\frac{1}{19}}), \quad a.s.\\
&\sup_{\y\in \cR^2} |n^{\frac{1}{2}}\beta_n(\y)-\bar{K}(TD_\y;n)|=O(n^{\frac{1}{2}\frac{2}{5}}),  
\quad a.s. 
\end{aligned}
\ee
where $\beta_n(\y)=n^{\frac{1}{2}}(F_n(\y)-F(\y))$ with 
$F_n(\y)$ the EDF and $D_\y=[0, y_1] \times [0, y_2]$.
\end{lemma}

\Rem{mks-2} We shall assume
$\sup_{\y}|\bar{B}_n(TD_\y)|=O(1)$ a.s. for the DF  $F$.
\end{remark}

   Here $\bar{B}$ and $\bar{K}$ are the stochastically equivalent versions of the ``measures" $B$ and $K$, and 
  R\'ev\'esz (1976)\cite{rp76} remarked that ``All the results here will be formulated and proved in the two-dimensional case only; it appears, however, that the generalization to higher dimensions is possible 
via the methods of this paper". 
To generalize the theorem to the d-dimensional case, we keep the definition of the Wiener Process $W(\x)=W(x_1,...,x_d)$ to be a separable Gaussian process from R\'ev\'esz's paper, and define 
$$
\begin{aligned}
&B.M.: \quad B(Q_z)=W(Q_z)-\lambda(Q_z)W(1,...,1)\\
&K.M.: \quad K(Q_z;y)=W(Q_z,y)-\lambda(Q_z)W(1,...,1,y)
\end{aligned}
$$  
where $\lambda(\cdot)$ is a Lebesgue measure on $\sZ^d$. The d-dimensional transformation is defined by Rosenblatt (1952)\cite{rm52}: Let $\X=(X_1, ..., X_d)$ be a random vector with DF $F(x_1,...,x_d)$. Let $\z=(z_1, ..., z_d)=T\x=T(x_1, ..., x_d)$, where $T$ is the transformation given by
$$
\begin{aligned}
& z_1=P(X_1\leq x_1)=F_1(x_1),\\
& z_2=P(X_2\leq x_2|X_1=x_1)=F_2(x_2|x_1),\\
& \vdots\\
& z_d=P(X_d\leq x_d|X_{d-1}\leq x_{d-1}, ..., X_1=x_1)=F_d(x_d|x_{d-1}, ..., x_1).
\end{aligned}
$$

We generalize the conditions to the d-dimensional case (S).
\begin{enumerate}
	\item[(S1)]
	$F(\x)$ is absolutely continuous on $\x\in \cR^d$.
	\item[(S2)] For all $\x=(x_1, ..., x_d) \in \cR^d$, there exists a constant $L>0$, 
	$$
	 \Big|\frac{\partial^2 F_2(x_2|F_1^{-1}(x_1))}{\partial x_1^2}\Big|\leq L, \quad
	\Big|\frac{\partial F_2(x_2|F_1^{-1}(x_1))}{\partial x_1}\Big|\leq L,
	$$
	\item[(S3)]  
	$$
	\begin{aligned}
	&\Big|\frac{\partial^2 F_3(x_3|F_2^{-1}(x_2|F_1^
		{-1}(x_1)),F_1^{-1}(x_1))}{\partial x_i \partial x_j}\Big|\leq L, \quad i,j=1,2,\\
	&\Big|\frac{\partial F_3(x_3|F_2^{-1}(x_2|F_1^
		{-1}(x_1)),F_1^{-1}(x_1))}{\partial x_i}\Big|\leq L, \quad i=1,2,
	\end{aligned}
	$$
	$$
	\vdots
	$$
	\item[(Sd)] For $d>2$,
	$$
	\begin{aligned}
	&\Big|\frac{\partial^2 F_d(x_d|F_{d-1}^{-1}(x_d|F_{d-2}^
		{-1}(x_{d-2}|...), ..., F_1^
		{-1}(x_1))}{\partial x_i \partial x_j}\Big|\leq L, \quad i,j=1, ..., d,\\
	&\Big|\frac{\partial F_d(x_d|F_{d-1}^{-1}(x_d|F_{d-2}^
		{-1}(x_{d-2}|...), ..., F_1^
		{-1}(x_1))}{\partial x_i }\Big|\leq L, \quad i=1, ..., d.
	\end{aligned}
	$$
\end{enumerate}
 
\Thm{drt} Let $\X_1=(X_{11},...,X_{1d})$, $\X_2=(X_{21},...,X_{2d})$, ... be a sequence of i.i.d. rv's having a common distribution function $F(\x)$. Assume (S). 
Then we can define a sequence $\{\bar{B}_n\}$ of Brownian Measures (B.M.) and a Kiefer Measure (K.M.) 
$\bar{K}$ such that almost surely, 
\bel{drt1}
\begin{aligned}
	&\sup_{\x\in \cR^d} |\beta_n(\x)-\bar{B}_n(TD_\x)|=O(n^{-\frac{1}{19}}), \\
	&\sup_{\x\in \cR^d} |n^{\frac{1}{2}}\beta_n(x)-\bar{K}(TD_\x;n)|=O(n^{\frac{1}{2}\frac{2}{5}}),  
\end{aligned}
\ee  
where $\beta_n(\x)=n^{\frac{1}{2}}(F_n(\x)-F(\x))$ with 
$F_n(\x)$ the EDF based on the sample $\X_1$, ..., $\X_n$, 
and $D_\x=[0, x_1]\times ... \times [0, x_d]$.

\end{theorem}

\def\Nu{{\mathcal N}}
\def\hbvt{\bdm{\check\btheta}}

\def\hbvt{\bdm{\hat\bvt}}
\def\ind{{\mathrm{ind}}}
 We need a property of U-statistics. 
 Let  $\bxi_1, \dots, \bxi_n$ be i.i.d. rv taking values in a measurable space $\cS$.  
 Let $\h$ be a measurable function from  $\cS^2$ to $\cR^m$ which is symmetric, i.e., 
$\h(\x, \y)=\h(\y, \x), \x, \y \in \cS$. A multivariate U-statistic (of order 2) with kernel $\h$ is 
defined as 
$$
\U_n(\h)={n \choose 2}^{-1}\sum_{1\leq i<j\leq n} \h(\bxi_i, \bxi_j).
$$
Assume that $\h$ is square-integrable. Let $\bmu(\h)=E(\h(\bxi_1, \bxi_2))$.  Recall that
a kernel $\k$ is \emph{degenerate} if $E(\k(\bxi_1, \bxi_2)|\bxi_2)=0$ a.s.   Let
$\bar \h(\x)=E(\h(\x, \bxi_2))$, and 
$$
\h^*(\x, \y)=\h(\x, \y)-\bar \h(\x)-\bar \h(\y)+\bmu(\h).
$$
Then $\h^*$ is a degenerate kernel. Let $\tilde \h=\h-\bmu(\h)$. Then 
$$
\h(\x, \y)=\bmu(\h)+\tilde \h(\x)+\tilde \h(\y)+\h^*(\x,\y).
$$
  One thus obtains the Hoeffding decomposition for a multivariate U-statistic, 
\bel{mhd}
\U_n(\h)=\bmu(\h)+\frac{2}{n}\sum_{j=1}^n \tilde \h(\bxi_j)+\U_n(\h^*)
=:\bmu(\h)+\hat\U_n(\h)+\U_n(\h^*), \quad a.s.
\ee
Let $\k$ be a degenerate kernel with $E(\|\k(\bxi_1, \bxi_2)\|^2)<\infty$. 
For $i<j, k<l$,  one has $E(\k(\bxi_i, \bxi_j)\k(\bxi_k, \bxi_l)^\top)=E(\k(\bxi_1, \bxi_2)^\vs)$
if $i=l, j=l$, and is equal to zero otherwise. Thus
\bel{mus-1}
E(\U_n(\k)^\vs)={n \choose 2}^{-1} E(\k(\bxi_1, \bxi_2)^\vs).
\ee
It is easy to see $E(\h^*(\bxi_1, \bxi_2)^\vs) \preceq E(\h(\bxi_1, \bxi_2)^\vs) $. Thus we prove 

\Lm{mus}  Suppose that $\h$ is a kernel with $E(\|\h(\bxi_1, \bxi_2)\|^2)<\infty$. Then 
$$
\U_n(\h)-\bmu(\h)-\hat\U_n(\h)=O_p(n^{-1} \sqrt{E(\|\h(\bxi_1, \bxi_2)\|^2)}).
$$
\end{lemma}

We need a Lipschitz-type property. 
\begin{itemize}
\item[(L)] Let $\tilde\bvt^{(i)}$ be  the estimator based on the observations with $\Z_i$ left out. 
Assume that there is a constant $L_0$ such that  
\bel{leave-out}
\max_{i}\|\tilde\bvt-\tilde\bvt^{(i)}\|\leq L_0/n. 
\ee
\end{itemize}
 Let $\tilde\bvt^{(ij)}$ be  the estimator based on the observations with $\Z_i ,\Z_j$ left out.
Applying (L) repeatedly, one has for some constant $L_0'$,  
\bel{leave-out2}
\max_{ij}\|\tilde\bvt-\tilde\bvt^{(ij)}\|\leq L_0'/n. 
\ee

Let $\v_n(\z_1)=\bPhi_{m_n}(F(\vep_1))\otimes \bPhi_{m_n} (G(\x_1))+
2(\h_{1, \A}(\z_1)+\h_{1, \B}(\z_1))$, where 
$$
\begin{aligned}
&\h_{1, \A}(\z_1)=E(\dot{\bPhi}_{m_n}( F(\vep_2))\otimes \bPhi_{m_n}(G(\X_2))(\1[\vep_1\leq \vep_2]-F(\vep_2)|\Z_1=\z_1),\\
&\h_{1, \B}(\z_1)=E({\bPhi}_{m_n}( F(\vep_2)) \otimes \dot{\Phi}_{m_n} (G(\X_2))(\1[\x_1\leq \X_2]-G(\X_2))).
\end{aligned}
$$

\Thm{mst} Suppose that the assumptions in \lmref{mst} hold. Assume (L), (S) and the assumptions
in \remref{resi} and \remref{mks-2}.  
Suppse that $\vep$ has  a bounded density.
Suppose that $\W_{n2}=E(\bPhi_{m_n}(G(\X))^\vs)$ is regular and that  
$\int \v\v^\top\,dQ$ is nonsingular.
If both $m_n$ and $n$ tend to infinity such that $m_n^{12}/n=o(1)$,  then $\hat\s_n$ satisfies the stochastic expansion,
\bel{mst6}
\hat \s_n=\s_n-\c\Var(\v(\Z))^{-1}\bar \v+o_p(n^{-1/2}),
\ee
where $\c=E\big(\w(\Z)\otimes \v^\top(\Z)\big)$. 
Thus, with $\D=\Var(\w(\Z))-\c \Var(\v(\Z))^{-1} \c^\top$, 
\bel{mst62}
\sqrt{n}(\hat\s_n-\bsigma(\bvt_0))\implies \Normal(0,\, \D),
\ee
As a consequence, $\hbvt$ given in \Ref{mst25} satisfies 
\bel{mst63}
\sqrt{n}(\hbvt-\bvt_0)\implies \Normal(0, \V),
\ee
where $\V=(\Delta_0^\top \H_0 \Delta_0)^{-1}\Delta_0^\top \H_0\D
\H_0 \Delta_0(\Delta_0^\top \H_0 \Delta_0)^{-1}$.

\end{theorem}

Proof of \thmref{mst}. We apply \thmref{4} with $\bpsi(\z)=\vecs((\z-\bmu)^\vs)$. 
Write $m=m_n$. 
As $\hat\u_n(\Z_j)=\hat{\g}(\Z_j)=\bPhi_m(F_n(\hat\vep))\otimes \bPhi_m (G_n(\X))$ and 
$m^7/n=o(1)$, 
\[
\maxj \|\u_n(\Z_j)\|+ \maxj \|\hat\u_n(\Z_j)\|\leq 4m^2=o(m^{-3/2}n^{1/2}).
\]
This shows \Ref{4a}. 
Since $\W_n=E(\u_n(\Z)\u_n(\Z)^\top)=\I_{m}\otimes \W_{n2}$ and $\bar{\W}_n=\jave \u_n(\Z_j)\u_n(\Z_j)^\top$, 
it follows that $\W_n$ is regular by the regularity of $\W_{n2}$, 
and that 
 $|\bar{\W}_n-\W_n|_o=O_p(m^2 n^{-1/2})$ as 
$$
\begin{aligned}
E|\bar{\W}_n-\W_n|_o^2
&\leq E\|\bar{\W}_n-\W_n\|^2
= \trace( E(\bar{\W}_n-\W_n)^\vs)\\
&\leq n^{-1} E\|\u_n(\Z)\|^4\leq m^4 n^{-1}.
\end{aligned}
$$
One verifies that there exists some constant $c_0>0$ such that for all $t$. 
 \bel{trigb}
\|{\bPhi}_{m}(t)\| \leq c_0 m^{1/2}, \quad
\|\dot{\bPhi}_{m}(t)\|\leq c_0 m^{3/2}, \quad
\|\ddot{\bPhi}_{m}(t)\|\leq c_0 m^{5/2}.
\ee
By the MVT, one thus has $|\hat{\W}_n-\bar\W|_o=O_p(m^5 n^{-1/2})$.
Taken together one proves $|\hat{\W}_n-\W_n|_o=o_p(m^{-1})$ as $m^{12}/n=o(1)$, yielding 
\Ref{4b}. Moreover, it is not difficult to verify that 
$\U_n=\W_n^{-1/2}\int \v_n\v_n^\top\,dQ\W_n^{-\top/2}=O(1)$.

  Write the left-hand-side average of \Ref{4c} as $\J_n+\K_n-E(\J_n+\K_n)$, where  
$$
\begin{aligned} 
& \J_n= \Ave \bpsi(\Z_j)\otimes(\hat\u_n(\Z_j)-\u_n(\Z_j)), \\
& \K_n= \Ave {\bpsi(\Z_j)\otimes (\u_n(\Z_j)-E(\u_n(\Z_j))}.
\end{aligned}
$$
Note first that 
\bel{mks-6}
E(\|\K_n\|^2) \leq n^{-1} E(|\bpsi(\Z)\|^2 \|\u_n(\Z)\|^2)=O(m^4n^{-1}).
\ee
We shall show next
\bel{mks-7}
E(\|\J_n\|^2)=O(m^4n^{-1}).
\ee
Taken together we prove \Ref{4c} as $m^5/n=o(1)$. 
To show \Ref{mks-7}, using the inequality $\|\A\otimes\B\|\leq \|\A\|\, \|\B\|$ and by \Ref{trigb},  we get 
\[
\begin{aligned}
\|\hat\u_n(\Z_j)-\u_n(\Z_j)\|
&\leq \|\bPhi_{m}(F_n(\hat\vep_j))-\bPhi_{m}(F(\vep_j))\|\cdot\|\bPhi_{m}(G_n(\X_j))\|\\
&\quad +\|\bPhi_{m}(F(\vep_j)\|\cdot\|\bPhi_{m}(G_n(\X_j))-\bPhi_{m}(G(\X_j))\|\\
&\leq c_0 m^2(|F_n(\hat\vep_j)-F(\vep_j)|+|G_n(\X_j)-G(\X_j)|).
\end{aligned}
\]
Let  $D_n=\sup_t |F_n(t)-F(t)|=O_p(n^{-1/2})$ (Kolmogorov-Simirnov's statistic). 
As $F$ has a bounded density (by $c_f$), we have 
\bel{ks}
|F_n(\hat\vep_j)-F(\vep_j)|
\leq D_n+|F(\hat\vep_j)-F(\vep_j)|
\leq D_n+c_f |\hat\vep_j-\vep_j|,
\ee
By \Ref{mks-3} below and \remref{resi}, we thus obtain 
\bel{mks-4}
 \Ave \|\hat\u(\Z_j)-\u(\Z_j)\|^2=O(m^4/n).
\ee
Therefore \Ref{mks-7} follows from 
\bel{mks-5}
E(\|\J_n\|^2) \leq E(\|\bpsi(\Z)\|^2)\Ave E(\|\hat\u(\Z_j)-\u(\Z_j)\|^2)=O(m^4/n).
\ee
  We shall now show \Ref{4d}--\Ref{4e}. Note 
\bel{a}
\begin{aligned}
&\bPhi_{m}(F_n(\hat\vep_j))=\bPhi_{m}(F(\vep_j))+\dot \bPhi_{m}(F(\vep_j))
(F_n(\hat\vep_j)-F(\vep_j))+\R_{1j}\\
&\bPhi_{m}(G_n(\X_j))=\bPhi_{m}(G(\X_j))+\dot \bPhi_{m}(G(\X_j))
(G_n(\X_j)-G(\X_j))+\R_{2j}, 
\end{aligned}
\ee
where,  by \Ref{ks} and the assumption in \remref{resi}, we have 
\bel{R1}
\maxj \|\R_{1j}\| 
=O_p(m^{5/2}/{n}).
\ee
By \Ref{mks} and the assumption in \remref{mks-2}, we have 
\bel{mks-3}
\maxj |G_n(\X_j)-G(\X_j)|=O_p(n^{-1/2}).
\ee
Similarly by \remref{mks-2},
\bel{R2}  
\max_j\|\R_{2j}\|=O_p(m^{5/2}/{n}).
 \ee
By \Ref{a}, 
\bel{118}
\Ave \hat \u_n(\Z_j)=\Ave \u_n(\Z_j)+\A+\B+\R,
\ee
where 
$$
\begin{aligned}
&\A=\Ave \dot{\bPhi}_m( F(\vep_j))\otimes \bPhi_{m}(G(\X_j))(F_n(\hat\vep_j)-F(\vep_j)), \\
&\B=\Ave {\bPhi}_m( F(\vep_j)) \otimes \dot{\Phi}_m (G(\X_j))(G_n(\X_j)-G(\X_j)),\\
&\R = \Ave \R_{1j} \otimes \bPhi_m(G_n(\X_j))+\Ave \bPhi_m(F_n(\hat\vep_j) \otimes \R_{2j}.
\end{aligned} 
$$
By \Ref{R1} and \Ref{R2} and the first equality in \Ref{trigb},  $\|\R\|=O(m^3/n)$.

Let $\b(\Z_i, \Z_j)={\bPhi}_m( F(\vep_j)) \otimes \dot{\Phi}_m (G(\X_j))(\1[\X_i\leq \X_j]-G(\X_j))$. 
It then follows $E(\kappa(\Z_i, \Z_j))=0$ for all $i, j$ from the independence of $\vep$ and $\X$,  and 
$\B$ is approximately a multivariate U-statistic, i.e., 
$\B=\U_n(\h_\B)+O(m^2/n)$, where $\h_\B(\z_1, \z_2)=\half (\b(\z_i, \z_j)+\b(\z_j, \z_i))$. 
Let $\h_1(\z_1)=E(\h(\z_1, \Z_2))$. Then
$$
\h_{1, \B}(\z_1)=E({\bPhi}_m( F(\vep_2)) \otimes \dot{\Phi}_m (G(\X_2))(\1[\x_1\leq \X_2]-G(\X_2))).
$$ 
By \lmref{mus}, 
\bel{decp1}
\B=\Ave 2\h_{1, \B}(\Z_j)+O_p(m^2/n).
\ee
Write $F_n(\hat\vep_j)-F(\vep_j)=(F_n(\hat\vep_j)-F_n(\vep_j))+(F_n(\vep_j)-F(\vep_j))$,  
giving $\A=\A_1+\A_2$.  Likewise, 
$\h_{1, \A}(\z_1)=E(\dot{\bPhi}_m( F(\vep_2))\otimes \bPhi_{m}(G(\X_2))(\1[\vep_1\leq \vep_2]-F(\vep_2)|\Z_1=\z_1)$, 
\bel{decp2}
\A_2=\Ave 2\h_{1, \A}(\Z_j)+O_p(m^2/n). 
\ee
    Using \Ref{leave-out}, one calculates $E(\|\A_{1}\|^2)=O(m^4/n^2)$. Hence 
\bel{29}
\A_{1}=\Ave \dot{\bPhi}_m( F(\vep_j))\otimes \bPhi_{m}(G(\X_j))(F_n(\hat\vep_j)-F_n(\vep_j))
=O_p(m^2/n).
\ee
This and \Ref{decp2} prove 
\bel{29-1}
\A=\Ave 2\h_{1, \A}(\Z_j)+O_p(m^2/n).
\ee
Taken together we show that \Ref{4e} holds with $\v=\u_n+2(\h_{1,\A}+\h_{1,\B})$ as $m_n^7/n=o(1)$. 
This also proves \Ref{4d}. \qed 

\section{Simulation results}
\label{sim}
We used the R package {\sc sem} to carry out the simulations based on the SEM \Ref{sems2}
with $\beta=1$.
The details of the package can be found in Fox (2006)\cite{fox06}. 
In LISREL notation and using Fig. 1, 
the SEM can be written as  
\bel{sem-22}
y_1 = \lambda_1x_2+\epsilon_1, \quad y_2=y_1+\lambda_3x_1+\lambda_2x_2+\epsilon_2.
\ee
The parameters to be estimated include all the regression coefficients $\lambda_1$, $\lambda_2$ and $\lambda_3$, and the measurement-error variances $\psi_1=Var(\epsilon_1)$ and $\psi_2=Var(\epsilon_2)$. We present the simulation results of estimating the coefficients $\lambda_1$, $\lambda_2$ and $\lambda_3$ with true parameter values $\lambda_1=1$, $\lambda_2=-1$ and $\lambda_3=0.5$. 

 For $n=30,50, 100$ and based on $50$ repetitions, we calculated the averages and medians of biases from
the usual and the EL-weighted MDF estimators $\lam_{n, k}$ and $\tilde\lam_k$ of the true parameter value $\lambda_k$, $k=1, 2, 3$, 
the averages and medians of the variances of the usual MDF and the EL-weighted estimators $v_{n,k}$ and $\tilde v_k$, $k=1, 2, 3$,  and the  
ratios $r_{1,k}=\bar{\tilde v}_k/\bar{v}_{n,k}$ and $r_{2,k}=median(\tilde v_k)/median(v_{n,k})$. 
The discrepancy function used is the ML discrepancy function given in \Ref{ml}.
A value of ratio less than one indicated the variance reduction of the EL-weighted
estimator over the usual estimator.  The results are reported on Tables \ref{t1}--\ref{t2}.  

 For Table \ref{t1}, the side information is the \emph{independence} of $\X$ and $\bep$ 
utilized via the constraint functions given in \Ref{sems4} for $m=1, 3, 5$,
where $\bep$ was generated from the normal mixture $0.9*N(0, \I_2)+0.1*N(0, 5\I_2)$, 
and $\X$ from the bivariate exponential $biexp(1, 3)$. 
 
 For Table \ref{t2}, the side information is \emph{known marginal medians} of $\X$ utilized via
the constraint functions given in \ref{sems5}, 
where $\X$ was generated from the bivariate exponential with scale parameters 
($\gamma_1, \gamma_2)$.  One can see that the efficiency gain is substantial (around 40\%). 
The ratios were stable with a slightly decreasing trend with increasing $n$, and the values of larger scale parameter had larger efficiency gains.

\begin{table}[!hbp]
	\centering
	\footnotesize\setlength{\tabcolsep}{2.5pt}
	\caption{\label{t1}
Simulated efficiency gain of the EL-weighted estimators in  the SEM \Ref{sems2}
using the side information in \Ref{sems4} of independence of $\varepsilon$ and $\X $ 
for a few values of $n$ and number $m$. 
$\bar{b}_{n,k}$ ($m(b_{n,k})$) and $\bar{\tilde{b}}_k$ ($m(\tilde{b}_k)$) are the averages (medians) of biases from
the usual and the EL-weighted MDF estimators respectively. $\bar{v}_{n,k}$ ($m(v_{n,k})$) and $\bar{\tilde{v}}_k$ ($m(\tilde{v}_k)$) are the averages (medians) of the variances of
the usual and the EL-weighted MDF estimators respectively. 
$r_{1,k}$ ($r_{2,k}$) are the ratios of the averages (medians) of the variances
of the EL-weighted estimators to the usual ones.  
$\bep\sim 0.9*N(0,\I_2)+0.1*N(0,5\I_2)$, $\X \sim \mathrm{biexp}(1,3)$. 
}
	\begin{tabular}{|c|c|c|c|c|c|c|c|c|c|c|c|}
		\hline
		\multicolumn{12}{|c|}{$n=30$}\\
		\hline
  & & & & & & & & & & & \\[-0.9em]
		$m$ & $\lambda_k$ & $\bar{b}_{n,k}$ & $\bar{\tilde{b}}_k$ & $m(b_{n,k})$ & $m(\tilde{b}_k)$ &  $\bar{v}_{n,k}$ & $\bar{\tilde{v}}_k$ &$r_{1,k}$& $m(v_{n,k})$ & $m(\tilde{v}_k)$& $r_{2,k}$\\
		\hline
		\multirow{3}{*}{$1$} 
		&	$\lambda_1$	&	0.2181	&	-0.2218	&	0.2876	&	-0.2437	&	2.4766	&	2.1962	&	0.8868	&	1.8529	&	1.3835	&	0.7467																
		\\
		&	$\lambda_2$	&	-0.0783	&	-0.0544	&	-0.1018	&	-0.0586	&	2.1048	&	1.9226	&	0.9134	&	1.6773	&	1.3444	&	0.8015																
		\\					
		&	$\lambda_3$	&	0.0333	&	-0.1002	&	-0.0855	&	-0.1314	&	0.5620	&	0.5660	&	1.0071	&	0.3729	&	0.3550	&	0.9520																
		\\
		\hline
		\multirow{3}{*}{$3$} 
		&	$\lambda_1$	&	-0.0893	&	-0.0688	&	0.0027	&	-0.0607	&	2.3394	&	1.7818	&	0.7616	&	1.5748	&	1.0300	&	0.6541																
		\\
		&	$\lambda_2$	&	0.1005	&	-0.1172	&	-0.3671	&	-0.3255	&	2.3664	&	1.8948	&	0.8007	&	1.8519	&	1.3011	&	0.7026																
		\\					
		&	$\lambda_3$	&	0.1303	&	-0.0938	&	0.0153	&	-0.0941	&	0.6303	&	0.5588	&	0.8866	&	0.5221	&	0.4538	&	0.8692																
		\\
		\hline
		\multirow{3}{*}{$5$} 
		&	$\lambda_1$	&	0.3742	&	-0.0233	&	0.2705	&	0.0381	&	1.7745	&	1.5109	&	0.8515	&	1.0320	&	0.8891	&	0.8615																
		\\
		&	$\lambda_2$	&	-0.1992	&	-0.1774	&	-0.1107	&	-0.1846	&	2.2048	&	1.5480	&	0.7021	&	1.0474	&	0.7961	&	0.7601																
		\\					
		&	$\lambda_3$	&	-0.0555	&	-0.1333	&	-0.0797	&	-0.0429	&	0.4901	&	0.3838	&	0.7831	&	0.2615	&	0.2143	&	0.8195																
		\\
		\hline
		\multicolumn{12}{|c|}{$n=50$}\\
		\hline
  & & & & & & & & & & & \\[-0.9em]
		$m$ & $\lambda_k$ & $\bar{b}_{n,k}$ & $\bar{\tilde{b}}_k$ & $m(b_{n,k})$ & $m(\tilde{b}_k)$ &  $\bar{v}_{n,k}$ & $\bar{\tilde{v}}_k$ &$r_{1,k}$& $m(v_{n,k})$ & $m(\tilde{v}_k)$& $r_{2,k}$\\
		\hline
		\multirow{3}{*}{$1$} 
		&	$\lambda_1$	&	0.0579	&	-0.0702	&	0.0963	&	-0.0546	&	1.4331	&	1.2840	&	0.8960	&	1.2614	&	0.9242	&	0.7327																
		
		\\
		&	$\lambda_2$	&	-0.1007	&	-0.1425	&	-0.1802	&	-0.1571	&	1.3882	&	1.2730	&	0.9170	&	1.1588	&	1.0470	&	0.9035																
		
		\\					
		&	$\lambda_3$	&	-0.2388	&	-0.1862	&	-0.1546	&	-0.1966	&	0.3278	&	0.3130	&	0.9549	&	0.2709	&	0.2552	&	0.9420																
		
		\\
		\hline
		\multirow{3}{*}{$3$} 
		&	$\lambda_1$	&	-0.0217	&	-0.0245	&	0.0697	&	-0.0398	&	1.2074	&	0.9621	&	0.7968	&	1.0255	&	0.7716	&	0.7524																
		
		\\
		&	$\lambda_2$	&	0.1368	&	-0.1279	&	0.1578	&	-0.1023	&	1.3264	&	1.0642	&	0.8023	&	1.0684	&	0.9074	&	0.8493																
		
		\\					
		&	$\lambda_3$	&	0.0039	&	-0.0187	&	-0.0181	&	-0.0118	&	0.2586	&	0.2127	&	0.8225	&	0.2367	&	0.1823	&	0.7702																
		
		\\
		\hline
		\multirow{3}{*}{$5$} 
		&	$\lambda_1$	&	0.1038	&	-0.0706	&	0.0412	&	-0.0716	&	1.1384	&	0.9720	&	0.8538	&	1.0185	&	0.8304	&	0.8153																
		
		\\
		&	$\lambda_2$	&	0.2014	&	-0.2178	&	0.1338	&	-0.1635	&	1.3568	&	1.1220	&	0.8269	&	1.1851	&	0.9699	&	0.8184																
		
		\\					
		&	$\lambda_3$	&	-0.0049	&	-0.0213	&	0.0461	&	-0.0224	&	0.3525	&	0.3362	&	0.9538	&	0.2649	&	0.2500	&	0.9438																
		
		\\
		\hline
		\multicolumn{12}{|c|}{$n=100$}\\
		\hline
  & & & & & & & & & & & \\[-0.9em]
		$m$ & $\lambda_k$ & $\bar{b}_{n,k}$ & $\bar{\tilde{b}}_k$ & $m(b_{n,k})$ & $m(\tilde{b}_k)$ &  $\bar{v}_{n,k}$ & $\bar{\tilde{v}}_k$ &$r_{1,k}$& $m(v_{n,k})$ & $m(\tilde{v}_k)$& $r_{2,k}$\\
		\hline
		\multirow{3}{*}{$1$} 
		&	$\lambda_1$	&	-0.0780	&	-0.0863	&	-0.0389	&	-0.0345	&	0.5978	&	0.5240	&	0.8765	&	0.5418	&	0.4167	&	0.7691																
		
		\\
		&	$\lambda_2$	&	0.0550	&	-0.0943	&	0.0390	&	-0.0564	&	0.6421	&	0.5687	&	0.8857	&	0.6456	&	0.4629	&	0.7170

		\\					
		&	$\lambda_3$	&	0.0102	&	-0.0581	&	0.0033	&	-0.0133	&	0.1508	&	0.1472	&	0.9761	&	0.1365	&	0.1281	&	0.9385

		\\
		\hline
		\multirow{3}{*}{$3$} 
		&	$\lambda_1$	&	0.1536	&	-0.0925	&	0.2299	&	-0.1106	&	0.5943	&	0.4759	&	0.8008	&	0.5237	&	0.4228	&	0.8073

		\\
		&	$\lambda_2$	&	-0.0769	&	-0.0637	&	-0.1581	&	-0.1517	&	0.6522	&	0.6210	&	0.9522	&	0.5397	&	0.5038	&	0.9335																
		
		\\					
		&	$\lambda_3$	&	0.0029	&	-0.0779	&	-0.0991	&	-0.0129	&	0.1621	&	0.1642	&	1.0130	&	0.1420	&	0.1410	&	0.9930																
		
		\\
		\hline
		\multirow{3}{*}{$5$} 
		&	$\lambda_1$	&	0.0736	&	-0.0769	&	0.1407	&	-0.1056	&	0.5228	&	0.4403	&	0.8422	&	0.4625	&	0.3472	&	0.7507																
		\\
		&	$\lambda_2$	&	-0.0601	&	-0.0712	&	-0.1616	&	-0.1752	&	0.5515	&	0.4704	&	0.8529	&	0.4872	&	0.3448	&	0.7077																
		\\					
		&	$\lambda_3$	&	-0.0090	&	-0.0242	&	-0.0339	&	-0.0770	&	0.1590	&	0.1465	&	0.9214	&	0.1299	&	0.1138	&	0.8761				
		\\
		\hline
	\end{tabular}
\end{table}

\begin{table}[!hbp]
	\centering
	\footnotesize\setlength{\tabcolsep}{2.5pt}
	\caption{\label{t2}
Same as Table \ref{t1} except for the side Information 
of known marginal medians of $\X$ with $\X$ generated from the bivariate exponential with scale 
parameters ($\gamma_1$, $\gamma_2$).}
	\begin{tabular}{|c|c|c|c|c|c|c|c|c|c|c|c|}
		\hline
		\multicolumn{12}{|c|}{$n=30, m=1$}\\
		\hline
  & & & & & & & & & & & \\[-0.9em]
		$(\gamma_1, \gamma_2)$ & $\lambda_k$ & $\bar{b}_{n,k}$ & $\bar{\tilde{b}}_k$ & $m(b_{n,k})$ & $m(\tilde{b}_k)$ &  $\bar{v}_{n,k}$ & $\bar{\tilde{v}}_k$ &$r_{1,k}$& $m(v_{n,k})$ & $m(\tilde{v}_k)$& $r_{2,k}$\\
		\hline
		\multirow{3}{*}{$(2,2)$} 
		&	$\lambda_1$	&	-0.0720	&	-0.0942	&	-0.1217	&	-0.1600	&	0.1731	&	0.1531	&	0.8845	&	0.1478	&	0.1333	&	0.9019															
		
		\\
		&	$\lambda_2$	&	0.0413	&	0.0298	&	-0.0104	&	0.0495	&	0.1934	&	0.1657	&	0.8568	&	0.1796	&	0.1425	&	0.7934															
		
		\\					
		&	$\lambda_3$	&	0.0453	&	0.0365	&	0.0885	&	0.0687	&	0.2067	&	0.1808	&	0.8747	&	0.1841	&	0.1640	&	0.8908
		\\	
		\hline
		\multirow{3}{*}{$(2,3)$} 
		&	$\lambda_1$	&	0.0424	&	0.0087	&	-0.0198	&	-0.0177	&	0.3424	&	0.2126	&	0.6209	&	0.3016	&	0.1918	&	0.6359															
		
		\\
		&	$\lambda_2$	&	-0.0937	&	-0.0527	&	-0.1107	&	-0.0952	&	0.4198	&	0.2852	&	0.6794	&	0.3813	&	0.2478	&	0.6499															
		
		\\					
		&	$\lambda_3$	&	-0.0600	&	-0.1512	&	-0.1224	&	-0.1159	&	0.2004	&	0.1819	&	0.9077	&	0.1626	&	0.1518	&	0.9336														
		
		\\
		\hline
		\multirow{3}{*}{$(2,4)$} 
		&	$\lambda_1$	&	0.0209	&	0.0535	&	-0.0327	&	-0.1286	&	0.7113	&	0.2411	&	0.3390	&	0.6064	&	0.2367	&	0.3903															
		
		\\
		&	$\lambda_2$	&	0.0860	&	0.1354	&	0.1475	&	0.1280	&	0.6681	&	0.3136	&	0.4694	&	0.5949	&	0.3068	&	0.5157														
		
		\\					
		&	$\lambda_3$	&	0.0825	&	0.0080	&	0.0064	&	-0.1096	&	0.1741	&	0.1917	&	1.1011	&	0.1312	&	0.1411	&	1.0755														
		
		\\
		\hline
		\multicolumn{12}{|c|}{$n=50$}\\
\hline
  & & & & & & & & & & & \\[-0.9em]
		$(\gamma_1, \gamma_2)$ & $\lambda_k$ & $\bar{b}_{n,k}$ & $\bar{\tilde{b}}_k$ & $m(b_{n,k})$ & $m(\tilde{b}_k)$ &  $\bar{v}_{n,k}$ & $\bar{\tilde{v}}_k$ &$r_{1,k}$& $m(v_{n,k})$ & $m(\tilde{v}_k)$& $r_{2,k}$\\
		\hline
		\multirow{3}{*}{$(2,2)$} 
		&	$\lambda_1$	&	0.0196	&	0.0318	&	0.0746	&	0.0918	&	0.0874	&	0.0721	&	0.8249	&	0.0797	&	0.0678	&	0.8507														
		
		\\
		&	$\lambda_2$	&	-0.0547	&	-0.0467	&	-0.0646	&	-0.0760	&	0.1021	&	0.0876	&	0.8580	&	0.0953	&	0.0830	&	0.8709

		\\					
		&	$\lambda_3$	&	-0.0847	&	-0.0766	&	-0.0894	&	-0.1072	&	0.0997	&	0.0876	&	0.8786	&	0.0877	&	0.0811	&	0.9247

		\\
		\hline
		\multirow{3}{*}{$(2,3)$} 
		&	$\lambda_1$	&	0.0059	&	0.0069	&	0.0058	&	0.0237	&	0.2121	&	0.1190	&	0.5611	&	0.1757	&	0.1072	&	0.6101															
		
		\\
		&	$\lambda_2$	&	-0.0544	&	-0.0326	&	-0.0447	&	-0.0061	&	0.2422	&	0.1489	&	0.6148	&	0.2288	&	0.1299	&	0.5677

		\\	&$\lambda_3$	&	0.0422	&	0.0342	&	0.0937	&	0.0515	&	0.1068	&	0.1026	&	0.9607	&	0.0954	&	0.0827	&	0.8669

		\\
		\hline
		\multirow{3}{*}{$(2,4)$} 
		&	$\lambda_1$	&	0.1288	&	0.1121	&	-0.0007	&	0.0403	&	0.3427	&	0.1501	&	0.4380	&	0.3068	&	0.1376	&	0.4485

		\\
		&	$\lambda_2$	&	0.1900	&	0.1183	&	0.0930	&	0.1413	&	0.3781	&	0.1822	&	0.4819	&	0.3453	&	0.1700	&	0.4923

		\\					
		&	$\lambda_3$	&	-0.0182	&	0.0027	&	-0.0527	&	0.0671	&	0.0978	&	0.0994	&	1.0164	&	0.0895	&	0.0867	&	0.9687	
		
		\\
		\hline
		\multicolumn{12}{|c|}{$n=100$}\\
\hline
  & & & & & & & & & & & \\[-0.9em]
		$(\gamma_1, \gamma_2)$ & $\lambda_k$ & $\bar{b}_{n,k}$ & $\bar{\tilde{b}}_k$ & $m(b_{n,k})$ & $m(\tilde{b}_k)$ &  $\bar{v}_{n,k}$ & $\bar{\tilde{v}}_k$ &$r_{1,k}$& $m(v_{n,k})$ & $m(\tilde{v}_k)$& $r_{2,k}$\\
		\hline
		\multirow{3}{*}{$(2,2)$} 
		&	$\lambda_1$	&	0.0353	&	0.0346	&	0.0460	&	0.0478	&	0.0369	&	0.0326	&	0.8837	&	0.0360	&	0.0326	&	0.9050	\\																								
		&	$\lambda_2$	&	-0.0511	&	-0.0585	&	-0.0298	&	-0.0316	&	0.0428	&	0.0378	&	0.8836	&	0.0423	&	0.0358	&	0.8456	\\																								
		&	$\lambda_3$	&	-0.0482	&	-0.0453	&	-0.0882	&	-0.0691	&	0.0431	&	0.0390	&	0.9045	&	0.0407	&	0.0361	&	0.8877	\\										\hline
		\multirow{3}{*}{$(2,3)$} 
		&	$\lambda_1$	&	-0.0809	&	-0.0784	&	-0.1066	&	-0.0858	&	0.0985	&	0.0677	&	0.6878	&	0.0941	&	0.0642	&	0.6819	\\																								
		&	$\lambda_2$	&	0.0191	&	0.0262	&	0.0175	&	0.0007	&	0.1037	&	0.0608	&	0.5865	&	0.1009	&	0.0673	&	0.6673	\\																								
		&	$\lambda_3$	&	-0.0104	&	-0.0128	&	-0.0179	&	-0.0125	&	0.0444	&	0.0436	&	0.9820	&	0.0395	&	0.0391	&	0.9899	\\										\hline
		\multirow{3}{*}{$(2,4)$} 
		&	$\lambda_1$	&	-0.0641	&	-0.0655	&	-0.0194	&	-0.0238	&	0.1732	&	0.1011	&	0.5838	&	0.1688	&	0.0973	&	0.5763	\\																								
		&	$\lambda_2$	&	0.0173	&	0.0186	&	-0.0162	&	0.0175	&	0.1747	&	0.1025	&	0.5868	&	0.1716	&	0.0848	&	0.4942	\\																								
		&	$\lambda_3$	&	-0.0484	&	-0.0464	&	-0.0168	&	-0.0152	&	0.0413	&	0.0409	&	0.9903	&	0.0385	&	0.0371	&	0.9636	\\										\hline

	\end{tabular}
\end{table}

\section{Proofs}
\label{proofs}

In this section, we first give two useful general theorems. As applications, we prove the
theorems presented in Section \ref{main}. 

Let $\x_1, ..., \x_n$ be $m$-dimensional vectors. Set
$$
\bar \x=\Ave \x_j, 
\quad x_*=\max_{1\leq j \leq n} \|\x_j\|, 
\quad \S=\Ave \x_j \x_j^\top,  
$$
$$
x^{(\nu)}=\sup_{\|\u\|=1} \Big\| \Ave (\u^{\top}\x_j)^{\nu}\Big\|,  \quad \nu=3,4,
$$
and let $\lambda$ and $\Lambda$ denote the smallest and largest eigen value 
of the matrix $\S$,
$$
\lambda = \inf_{\|\u\|=1} \u^{\top}\S \u  
\und
\Lambda = \sup_{\|\u\|=1} \u^{\top} \S \u.
$$
With these we associate the empirical likelihood
$$
\mathscr{R} = \sel{ \pave \x_j=0}. 
$$

\PS\ \cite{PS12} carefully examined the above maximization as a numeric problem and detailed some very useful properties. We quote their Lemma 5.2 below for our application.

\Lm{owen}
The inequality $\lambda > 5\|\bar \x\| x_*$ implies that there is a unique 
$\bzeta$ in $\cR^m$ satisfying the below \Ref{o1} to \Ref{o9}.

\bel{o1}
1+\bzeta^{\top}\x_j>0, \quad j=1,\dots,n,
\ee
\bel{o2}
\sum_{j=1}^n \dfrac{\x_j}{1+\bzeta^{\top} \x_j}=0,
\ee

\bel{o3}
\|\bzeta\|\leq \frac{\|\bar \x\|}{\lambda -\|\bar \x\| x_*},
\ee

\bel{o4}
\|\bzeta\|x_*\leq \frac{\|\bar \x\|x_*}{\lambda -\|\bar \x\| x_*} < \frac{1}{4},
\ee
%
%
\bel{o6}
\Ave (\bzeta^{\top}\x_j)^2 = \bzeta^{\top} \S \bzeta 
\le \Lambda \|\bzeta\|^2 
\leq \dfrac{\Lambda\|\bar \x\|^2}{(\lambda- \|\bar \x\|x_*)^2},
\ee
%
\bel{o8}
\Big\|\Ave \dfrac{\r_i}{1+\bzeta^{\top}\x_j}-\r_j+\r_j \x_j^\top \bzeta \Big\|
\leq \Big\| \Ave \r_j (\bzeta^{\top}\x_j)^2 \Big\| + 
\dfrac{4}{3} \Ave \|\r_j\|\; |\bzeta^{\top} \x_j|^3,
\ee
for vectors $\r_1,\dots,\r_n$ of the same dimension, and 
\bel{o9}
\|\bzeta-\S^{-1}\bar \x\|^2 \leq  2 \Big(\dfrac{1}{\lambda} + \dfrac{\Lambda}{9\lambda^2} \Big) 
\|\zeta\|^4 x^{(4)}.
\ee

\end{lemma}

Now use the fact that $\|\x\|=\sup_{\|\v\|=1} \v^{\top}\x$, 
the Cauchy-Schwartz inequality,
\Ref{o3}, \Ref{o4} and \Ref{o6} to bound the square of the first term 
of the right-hand side of \Ref{o8} by 
$$
\Ave (\bzeta^{\top} \x_j)^4 \sup_{\|\v\|=1} \v^{\top} \Big(\Ave \r_j\r_j^\top\Big) \v  
\le  \|\bzeta\|^4 x^{(4)}\Big|\Ave \r_j\r_j^\top\Big|_o
$$
and the square of the second term by 
$$
\frac{16}{9} x_*^2  \|\bzeta\|^2 \Ave \|\r_j\|^2
\Ave (\bzeta^{\top} \x_j)^4
\le \frac{16d}{9} (x_* \|\bzeta\|)^2 \|\bzeta\|^4 \x^{(4)}\Big|\Ave \r_j\r_j^\top\Big|_o,
$$
where $d$ is the dimension of $r_j$. 
Combining the above we obtain 
\bel{o82}
\Big\|\Ave\! \dfrac{\r_i}{1+\bzeta^{\top}\x_j}-\r_j+\r_j \x_j^\top \bzeta \Big\|^2\!\!
\!\leq \!\|\bzeta\|^4 x^{(4)}\!\Big|\Ave\! \r_j\r_j^\top\!\Big|_o\big[1+\frac{16d}{9} (x_* \|\bzeta\|)^2\big].
\ee

We now apply the above results to random vectors.
Let $\T_{n1},\dots, \T_{nn}$ be $m_n$-dimensional random vectors. 
With these random vectors we associate the empirical likelihood
$$
\el =\sel{\pave \T_{nj}=0}.
$$
To study the asymptotic behavior of $\el$ we introduce 
$$
T_n^*= \maxj \|\T_{nj}\|,  \quad \bar \T_n = \ave \T_{nj}, \quad T_n^{(\nu)}=\sup_{\|\u\|=1}
\Big\|\Ave (\u^{\top} \T_{nj})^{\nu}\Big\|, 
$$
and the matrix
$$
\S_n = \Ave \T_{nj}\T_{nj}^{\top},
$$
and let $\lambda_n$ and $\Lambda_n$ denote the smallest 
and largest eigen values of $\S_n$,
$$
\lambda_n = \inf_{\|\u\|=1} \u^{\top}\S_n \u  
\und
\Lambda_n = \sup_{\|\u\|=1} \u^{\top}\S_n \u.
$$

We impose the following conditions on $\T_{nj}$.
\begin{enumerate}
	\item[(A1)]
	$T_n^*= o_p(m_n^{-3/2}n^{1/2})$.
	\item[(A2)]
	$\|\bar \T_n\| = O_p(m_n^{1/2}n^{-1/2})$.
	\item[(A3)] 
	There is a sequence of regular $m_n\times m_n$ dispersion matrices $\W_n$  such that 
	$$
	|\S_n - \W_{n}|_o = o_p(m_n^{-1}).
	$$
\end{enumerate}

We impose the following conditions on $\bpsi$ and $\T_{nj}$.

\begin{enumerate}
	\item[(B1)]$\ave\big(\bpsi(Z_j)\otimes \T_{nj}^\top-E\big(\bpsi(Z_j)\otimes \T_{nj}^\top\big)\big)=o_p(m_n^{-1/2})$.
	\item[(B2)] There exists some measurable function $\bchi$ from $\cZ$ into $\cR^{d}$ such that 
	$\int \bchi \,dQ=0$, $\int \|\bchi\|^2\,dQ<\infty$ and
	$$
	\iAve\pare{ \A_n\W_n^{-1}\T_{ni}-\bchi(Z_i)}=o_p(n^{-1/2}),
	$$
	where $\A_n=:\ave E\big(\bpsi(Z_j)\otimes \T_{nj}^\top\big)$.
\end{enumerate}
Let us first consider the case that $m_n$ tends to infinity with the sample size.  
We have the following result.
\Thm{5}
Suppose (A1)-(A3) and (B1)-(B2) hold. 
Then there exists unique $\bzeta_n$ which satisfies
\bel{5.1}
1+\bzeta_n^\top \T_{nj}>0, \quad
\Ave \frac{\T_{nj}}{1+\bzeta_n^\top \T_{nj}}=0,
\end{equation}
such that as $m_n$ tends to infinity, 
\bel{5.2}
\btheta_n=:\Ave \frac{\bpsi(Z_j)}{1+\bzeta_n^\top \T_{nj}}
=\bar \bpsi-\bar \bchi + o_p(n^{-1/2}), 
\end{equation}
where $\bar \bchi=\ave \bchi(Z_j)$ with $\bchi$ given in (B2). 
\end{theorem}
\proof It follows from (A1) and (A2) that $T_n^*\|\bar \T_n\| = o_p(1)$, 
and from (A3) that there are positive numbers $a<b$ such that 
$P(a \leq \lambda_n \leq \Lambda_n \leq b) \to 1$.
Thus all three conditions imply that the probability of the event 
$\{\lambda_n > 5 T_n^* \|\bar \T_n\|\}$ tends to one. 
Consequently, by \lmref{owen},  there exists an $m_n$-dimensional random vector $\bzeta$ 
which is uniquely determined on this event by the properties \Ref{o1}--\Ref{o82} including \Ref{5.1}. 
To prove \Ref{5.2}, we apply \Ref{o82} with $\r_j=\bpsi(Z_j)$. 
Note first that 
\bel{5.3}
T_n^{(4)} \le \Lambda_n (T_n^*)^2.
\end{equation}
This, \Ref{o3} and (A1)-(A2) imply that the right side of \Ref{o82} is bounded by  
$$
\Big(1+\frac{d}{9}\Big)\frac{\|\bar \T_n\|^4}{(\lambda_n-\|\bar \T_n\|T_n^*)^4} \Lambda_n(T_n^*)^2
\big|\Ave \bpsi(Z_j)\bpsi(Z_j)^\top\big|_o=o_p(m_n^{-1}n^{-1}),
$$
where the equality holds since the spectral norm of the average is bounded due to the square-integrability of $\bpsi$. Thus from \Ref{o82} it follows  
\bel{5a}
\btheta_n=\Ave \bpsi(Z_j)-\Ave \bpsi(Z_j)\otimes \T_{nj}^\top \bzeta_n +o_p(n^{-1/2}).
\end{equation}
In view of (B2) the desired \Ref{5.2} now follows from  \Ref{5b}-\Ref{5d} below.
\bel{5b}
\Ave \pare{\bpsi(Z_j)\otimes \T_{nj}^\top -E\big(\bpsi(Z_j)\otimes \T_{nj}^\top\big)}\bzeta_n=o_p(n^{-1/2}),
\end{equation}

\bel{5c}
\Ave E\big(\bpsi(Z_j)\otimes \T_{nj}^\top\big)\big(\bzeta_n-\S_n^{-1}\bar \T_n\big)=o_p(n^{-1/2}),
\end{equation}

\bel{5d}
\Ave E\big(\bpsi(Z_j)\otimes \T_{nj}^\top\big)\big(\S_n^{-1}-\W_n^{-1}\big)\bar \T_n=o_p(n^{-1/2}).
\end{equation}
Note first that (B1), (A2) and \Ref{o3} imply \Ref{5b}. Next we show
$$
\A_n=\Ave E\big(\bpsi(Z_j)\otimes \T_{nj}^\top\big)=O(m_n^{1/2}).
$$ 
Indeed, by Cauchy inequality, 
\begin{eqnarray*}
	&& \|\A_n\|^2\le \Ave \|E\big(\bpsi(Z_j)\otimes \T_{nj}^\top\big)\|^2
	\le E\big(\|\bpsi(Z_1)\|^2\big)\Ave E(\|\T_{nj}\|^2\big)\\
	&&\hspace{5.5cm} =E\big(\|\bpsi(Z_1)\|^2\big)\trace\pare{E(\S_n)}.
\end{eqnarray*}
But by (A3), the above trace is bounded by 
$$
\|\trace\pare{E(\S_n-\W_n)}\|+\trace\pare{E(\W_n)}
\le m_n E\pare{|\S_n-\W_n|_o}+\Lambda_n m_n,
$$
Thus $\|\A_n\|^2=O(m_n)$.  This, the regularity of $\W_n$ in (A3), \Ref{o9},  \Ref{5.3} and (A1) imply that 
the square of the right side of \Ref{5c} is bounded by
$$
\begin{aligned}
O(m_n) O_p(\|\bzeta_n\|^4 T_n^{(4)})
&=O(m_n) O_p(\|\bar \T_n\|^4 (T_n^{*})^2)\\
&=O(m_n) o_p(m_n^{2}n^{-2} m_n^{-3} n)
=o_p(m_n^{-1}n^{-1}),
\end{aligned}
$$ 
hence \Ref{5c} is proved. Again the rate of $A_n$ and (A2)-(A3) imply \Ref{5d}. 
This completes the proof.
\endproof

Examining the proof of \thmref{5} one can see  the following holds. 
\Thm{6} Suppose (A1)-(A3) and (B1)-(B2) are met for fixed $m_n=m$.
Then the results in \thmref{5} hold as $n$ tends to infinity. 
\end{theorem}

\medskip  
{\sc Proof} of \thmref{1}. We verify the conditions of \thmref{6} with $\T_{nj}
=\u(Z_j)$. Since $\u$ is square-integrable, conditions (A1) -- (A3) are  satisfied
with $\W_n=\W=E(\u\u^\top(Z))$. The square-integrability of $\bpsi$ implies that 
(B1) -- (B2) are met with $\A_n=\A=E(\bpsi(Z)\otimes \u(Z)^\top)$ and 
$\bchi=\A\W^{-1}\u$, by the weak law of large numbers. We now apply the result 
of \thmref{6} to complete the proof. 
\qed

\medskip

{\sc Proof of \thmref{2}.} We shall apply \thmref{6} to prove the result with $\T_{nj}=\hat \u(Z_j)$.
Clearly conditions (A1), (A3) and (B1) follows from \Ref{2a} -- \Ref{2c} respectively,
whereas (A2) is implied by \Ref{2e} in view of the fact that the right-hand-side average
of \Ref{2e} is $O_p(n^{-1/2})$.   
By Cauchy inequality,
$$
\big\|\Ave E\big(\bpsi(Z_j)\otimes (\hat \u(Z_j)-\v(Z_j))\big)\big\|^2
\leq E(\|\bpsi(Z_1)\|^2) \Ave E\big(\|\hat \u(Z_j)-\v(Z_j)\|^2\big),
$$
which is $o(1)$ by \Ref{2d}. Hence the $\A_n$ in (B2) is given by 
$$
\A_n=\Ave E\big(\bpsi(Z_j)\otimes \hat \u(Z_j)\big)
=E(\bpsi(Z_1)\v(Z_1))+o(1). 
$$
Thus by \Ref{2e}, (B2) holds with $\bchi=E(\bpsi(Z_1)\v(Z_1))\W^{-1}\v$. 
We now apply \thmref{6} to complete the proof.
\qed

\medskip 
{\sc Proof} of \thmref{3}. 
We apply \thmref{5} with $\T_{nj}=\u_n(Z_j)$ to prove the result, i.e. 
verify its conditions (A1)-(A3) and (B1)-(B2). Obviously 
\Ref{3a}, \Ref{3b} and \Ref{3c} correspond to (A1), (A3) and (B1) respectively. 
It follows from the regularity of $\W_n$ that $\trace(\W_n)\le B m_n$ for some constant $B$. Thus from $nE[\|\bar \T_n\|^2] = \trace(\W_n)=O(m_n)$ it yields (A2). We are now left to prove (B2). 
Noticing $\A_n=E(\bpsi(Z)\otimes\u_n^\top(Z))$ and $\W_n=E(\u_n\u_n(Z)^\top)$, and
$\A_n\W_n^{-1}\u_n$ is the projection of $\bpsi(Z)$ onto the closed linear span $[\u_n]$, so 
that $\A_n\W_n^{-1}\u_n(Z)$ is the conditional expectation of $\bpsi(Z)$ given $\u_n(Z)$, i.e., 
$$
\A_n\W_n^{-1}\u_n(Z)=E(\bpsi(Z)|\u_n(Z)).
$$
Since $E(\bpsi(Z)|\u_n(Z)), n\ge 1 $ forms a martingale with respect to the sigma algebras,   
$\sigma(\u_n(Z)), n\ge 1$, generated by $\u_n(Z)$, it follows from L\'evy's martingale convergence theorem (see e.g. page 510, Shiryaev \cite{sh96}) that 
$$
E(\bpsi(Z)|\sigma(\u_n(Z))) \to E(\bpsi(Z)|\sigma(\u_\infty(Z))), \quad a.s. \quad n\to\infty.
$$
By the property of conditional expectation (see e.g. Proposition 1, page 430, Bickel, {\sl et al.} 
\cite{bk93}), the last conditional expectation is the projection
of $\bpsi(Z)$ onto the closed linear span $[\u_\infty(Z)]$, i.e.,
$E(\bpsi(Z)|\sigma(\u_\infty(Z)))=\Pi(\bpsi(Z)|[\u_\infty(Z)])$, hence
$$
\bvarphi(Z)=\Pi(\bpsi(Z)|[\u_\infty(Z)])=E(\bpsi(Z)|\u_\infty(Z)).
$$
Thus that (B2) is satisfied with $\bchi=\bvarphi$ follows from  
$$
nE\Big(\|\iAve\pare{\A_n\W_n^{-1}\u_n(Z_i)-\bvarphi(Z_i)}\|^2\Big)
=E\Big(\|\A_n\W_n^{-1}\u_n(Z)-\bchi(Z)\|^2\Big), 
$$
which converges to zero as $n$ tends to infinity by the property of the convergence of
Fourier series. This completes the proof.
\qed

\medskip
{\sc Proof} of \thmref{4}. We prove the result by verifying conditions (A1)-(A3) and (B1)-(B2) of \thmref{5} with $\T_{nj}=\hat\u_n(Z_j)$. Clearly $(A1)$, $(A3)$ and (B1) correspond to \Ref{4a}, \Ref{4b} and \Ref{4c} respectively, while (A2) follows from \Ref{4d} and \Ref{pf41} below.  
We are left to verify (B2). First, by Cauchy inequality and \Ref{4d}, 
\begin{eqnarray*}
	&& \Big\|\Ave E\pare{\bpsi(Z_j)\otimes (\hat\u_n(Z_j)-\v_n(Z_j)}\Big\|^2\\
	&& \leq E(\|\bpsi(Z)\|^2) \Ave E\pare{\|\hat\u_n(Z_j)-\v_n(Z_j)\|^2}=o(m_n^{-1}), 
\end{eqnarray*}
so that the $\A_n$ in (B2) satisfies
$$
\A_n=E\pare{\bpsi(Z)\otimes \v_n(Z)}+o(m_n^{-1/2}).
$$
Note that $\trace(\U_n)\le m_n |\U_n|_o= O(m_n)$ and  
$$
\begin{aligned}
nE(\|\bar\v_n\|^2)
&=E(\|\v_n(Z)\|^2)\leq |\W_n^{1/2}|_o^2 E(\|\W_n^{-1/2}\v_n(Z)\|^2)\\
& = |\W_n^{1/2}|^2_o \trace(\U_n).
\end{aligned} 
$$
This shows 
\bel{pf41}
\|\bar \v_n\|= O_p(n^{-1/2}m_n^{1/2}).
\ee 
By Cauchy inequality, 
$$
\|E\pare{\bpsi(Z)\otimes \v_n(Z)}\|^2 \leq 
E\big(\|\bpsi(Z)\|^2\big)E\big(\|\v_n(Z)\|^2\big).
$$
But 
$$
\begin{aligned}
E\big(\|\v_n(Z)\|^2\big)
&\leq |\W_n^{1/2}|_o^2E\big(\|\W_n^{-1/2}\v_n(Z)\|^2\big)\\
&=|\W_n^{1/2}|_o^2\trace(\U_n)
=O(m_n).
\end{aligned}
$$
Hence 
$$
E\pare{\bpsi(Z)\otimes \v_n(Z)}=O(m_n^{1/2}).
$$
Therefore combining the above and in view of \Ref{4e} we arrive at  
$$
\begin{aligned}
\Ave \A_n\W_n^{-1}\hat \u_n(Z_j)
&=\big(E\pare{\bpsi(Z)\otimes \v_n(Z)}+o(m_n^{-1/2})\big)\W_n^{-1}\\
& \hspace{2.7cm} \times \big(\bar \v_n+o_p(m_n^{-1/2}n^{-1/2})\big)\\
&= E\pare{\bpsi(Z)\otimes \v_n(Z)}\W_n^{-1}\bar\v_n +o_p(n^{-1/2}),
\end{aligned}
$$
Analogous to the proof of (B2) in \thmref{3} (or applying $\u_n=\v_n$), we have
$$
E\pare{\bpsi(Z)\otimes \v_n(Z)}\W_n^{-1}\bar\v_n=\bar \bchi+o_p(n^{-1/2}),
$$ 
where $\bchi=\Pi(\bpsi|[\v_\infty])$ is the projection of $\bpsi(Z)$ onto the closed linear
span $[\v_\infty]$. Clearly $\int \bchi\,dQ=0$ and $\int |\bchi|^2\,dQ<\infty$. 
Thus (B2) is proved with $\bvarphi=\bchi=\Pi(\bpsi|[\v_\infty])$. 
This finishes the proof.
\qed

\end{document}